\documentclass{amsart}

\usepackage{eufrak}

\usepackage{amssymb,amsmath, amscd, mathrsfs, yfonts, bbm, graphics, dsfont}

\title[Non-Commutative Functions, Free Levy-Hincin Formula] {Non-Commutative Functions and Non-Commutative Free Levy-Hincin Formula}

\author{Mihai Popa and Victor Vinnikov}

\address{Center for Advanced Studies in Mathematics at the Ben Gurion University
of Negev, P.O. B. 653, Be�er Sheva 84105, Israel and
\newline
Institute of Mathematics �Simion Stoilow� of the Romanian Academy, P.O. Box 1-764,
Bucharest, RO-70700, Romania}
\email{popa@math.bgu.ac.il}
\address{Department of Mathematics, Ben Gurion University of Negev, Be�er Sheva
84105, Israel}
\email{vinnikov@cs.bgu.ac.il}

\newtheorem{claim}{}[section]
\newtheorem{defn}[claim]{Definition}
\newtheorem{thm}[claim]{Theorem}
\newtheorem{lemma}[claim]{Lemma}
\newtheorem{remark}[claim]{Remark}
\newtheorem{prop}[claim]{Proposition}
\newtheorem{cor}[claim]{Corollary}

\newcommand{\cA}{\mathcal{A}}
\newcommand{\cH}{\mathcal{H}}

\newcommand{\cT}{\mathcal{T}}
\newcommand{\cB}{\mathcal{B}}
\newcommand{\cD}{\mathcal{D}}
\newcommand{\cC}{\mathcal{C}}
\newcommand{\cR}{\mathcal{R}}
\newcommand{\Sb}{\Sigma_{\mathcal{B}}}
\newcommand{\Sbd}{\Sigma_{{\mathcal{B}:\mathcal{D}}}}
\newcommand{\cK}{\mathcal{K}}

\newcommand{\id}{\text{Id}}

\newcommand{\cfree}{\framebox{c}}
\newcommand{\cka}{{}^c\kappa}
\newcommand{\vu}{\vec{u}}
\newcommand{\CR}{{}^cR}

\newcommand{\cE}{\mathcal{E}}

\newcommand{\cL}{\mathcal{L}}

\newcommand{\lra}{\longrightarrow}

\newcommand{\X}{\mathcal{X}}

\newcommand{\bx}{\cB\langle \mathcal{X} \rangle}

\usepackage{amssymb}
\input xy
\xyoption{all}

\newcommand{\ncspace}[1]{\ensuremath{#1}_{\text{nc}}}
\newcommand{\mat}[2]{\ensuremath{{#1}^{#2\times #2}}}
\newcommand{\rmat}[3]{\ensuremath{{#1}}^{#2\times #3}}
\newcommand{\ball}{\mathbb{B}}
\newcommand{\Nilp}{\text{Nilp}}
\newcommand{\ten}[1]{\mathbf{T}(#1)}
\newcommand{\talpha}{\widetilde{\alpha}}
\newcommand{\tphi}{\widetilde{\phi}}
\newcommand{\tmu}{\widetilde{\mu}}

\newcommand{\cV}{\mathcal{V}}
\newcommand{\cW}{\mathcal{W}}
\begin{document}

\maketitle
\begin{abstract}The paper is discussing infinite divisibility in the setting of operator-valued boolean, free and, more general, c-free independences. Particularly, using Hilbert bimodules and non-commutative functions techniques, we obtain analogues of the Levy-Hincin integral representation for infinitely divisible real measures.
\end{abstract}

\section{Introduction and notations}

The paper presents some results concerning infinite divisibility in the framework of operator-valued non-commutative probability.

In probability theory - classic and non-commutative - limit theorems play a central role. The ``most general'' limit theorems involve so-called infinitesimal arrays, and the limit distributions are usually identified with {\em infinitely divisible} distributions.
There is a consistent literature about infinitely divisible measures in classical probability  (see \cite{feller}), dating back to Kolmogorov (\cite{kolmogorov}), P. Levy (\cite{levy}) etc. Similar results have been found for non-commutative independences, such as \cite{BP} and \cite{BV} for free independence, \cite{speicherwaroudi} for boolean probability and \cite{Krystek}, \cite{Wang} for conditionally free probability. In the operator-valued case, when states are replaced by positive conditional expectations or, more general, by completely positive maps between C$^\ast$- or operator algebras, very little was known about operator-valued infinite divisibility; the only exception we know of was Speicher's work \cite{speicherhab}. One of the obstructions is that while in the scalar case important results characterizing infinite divisibility are coming from Nevalinna-Pick representation properties of the functions that linearize additive convolutions (such as the log of the Fourier transform in the classic case or the Voiculescu's $R$- and $\phi$-transforms in the free case), such analytic tools are not yet available in the operator-valued case. The new topic of non-commutative (\cite{ncfound}) or completely matricial (\cite{voi09}) functions may possibly fill this gap. In particular, the results from Section 5 of the present paper indicate that results similar to the Nevalinna-Pick representation hold also for non-commutative functions.

The paper is organized in five sections. First section presents the introduction and some notations. Sections 2-4 are aimed towards results characterizing infinite divisibility in operator-valued non-commutative probability using combinatorial and operator-algebras methods and constructions. More precisely, in Section 2 we use a non-commutative version of the ``Boolean Fock space'' construction from \cite{anshelevich} to prove that, as in the scalar case, boolean infinite divisibility is trivial in the operator-valued case; particularly, any completely positive map between two C$^\ast$-algebras is boolean infinitely divisible. Section 3 is describing infinite divisibility with respect to free independence over a positive conditional expectation in terms of maps satisfying a condition of complete positivity. Section 4 is utilizing the techniques from the Boolean case (Section 2) to extend the results of Section 3 from positive conditional expectations to completely positive maps. In particular, we present a construction of the non-commutative version of the conditionally free $R$-transform of Bozejko, Leinert and Speicher in terms of creation, annihilation and preservation operators on a certain inner-product bimodules. In the scalar-valued case this construction gives a new, combinatorial proof of the main result from \cite{Krystek} (see also \cite{Wang}) characterizing conditionally free infinite divisibility. In Section 5 we use the tools from the theory of non-commutative functions (see \cite{ncfound}, \cite{voi09}) to define the non-commutative $R$- (also constructed in \cite{voi09}), ${}^cR$- and $B$-transforms. Reformulated in terms of these transforms, the results from Sections 2, 3 and 4 are very similar to the free and conditionally free versions of the Levy-Hincin formula from \cite{BP}, respectively \cite{Krystek} and \cite{Wang}. The present material is using the notions detailed in \cite{ncfound} (work in progress), but it is self-contained in this regard, the needed material on non-commutative functions is briefly discussed in Section 5.2.

 We will introduce now several notations. Throughout the paper $\cB$ will be a unital C$^\ast$-algebra. We will denote by $\bx$  the $\ast$-algebra freely generated by $\cB$ and the selfadjoint symbol $\X$. Unless otherwise explicitly stated, we do not suppose that $\cB$ commutes with $\X$.  We will also use the notations $\cB\langle  X \rangle _0$ for the $\ast$-subalgebra of $\bx$ of all polynomials without a free term, and the notation $\cB\langle \X_1, \X_2, \dots\rangle$ for the $\ast$-algebra freely generated by $\cB$ and the non-commutating self-adjoint symbols $\X_1, \X_2, \dots$.

In several instances we will identify $\mathcal{T}(\cB)$, the tensor algebra over $\cB$, to the subalgebra of $\bx$ spanned by $\{\X b_1\X b_2\cdots \X b_n\colon  n\in\mathbb{N}, b_1,\dots,b_n\in\cB\}$  via
\[b_1\otimes b_2\otimes\cdots\otimes b_n\mapsto \X b_1\X b_2\cdots \X b_n.\]

 The set of all positive conditional expectations from $\bx$ to $\cB$ will be denoted by $\Sb$.

 For $\cB\subseteq\cD$ a unital inclusion of C$^\ast$-algebras, we denote by $\Sbd$ the set of all unital,  positive $\cB$-bimodule maps $\mu:\bx\lra\cD$ with the property that for all positive integers $n$ and all $\{f_i(\X)\}_{i=1}^n\subset \bx$ we have that:
 \begin{equation}\label{prop} \bigl[\mu(f_j(\X)^\ast f_i(\X))]_{i,j=1}^n\geq 0\ \text{in $M_n(\cD)$}.
 \end{equation}

 Remark that $\Sb=\Sigma_{\cB:\cB}$, as an easy consequence of Exercise 3.18 from \cite{paulsen}.

 We will denote by $\Sb^0$, respectively $\Sbd^0$, the set of all $\mu\in\Sb$ (respectively $\mu\in\Sbd$ whose moments do not grow faster than exponentially, that is there exists some $M>0$ such that for all $b_1, \dots, b_n\in \cB$, we have
\begin{equation}\label{pr2} ||\mu(\X b_1 \X b_2\cdots \X b_n \X)||< M^{n+1}||b_1||\cdots ||b_n||.
 \end{equation}

  We will also use the following definition (see \cite{lance}):

 \begin{defn}
  Let $A$ be a C$^\ast$-algebra. A \emph{semi-inner-product $A$-module} is a linear space $E$ which is a right $A$-module together with a map $(x,y)\mapsto \langle x, y\rangle:E\times E\lra A$ such that ($x, y, z\in E, a\in A, \alpha, \beta\in\mathbb{C}$):
  \begin{enumerate}
  \item[(i)]\ $\langle\alpha x+\beta y, z \rangle=\alpha\langle x, y\rangle+\beta\langle y, z\rangle$
 \item[(ii)]\ $\langle xa , y\rangle=\langle x , y\rangle a$
  \item[(iii)]\ $\langle x, y\rangle=\langle y, x\rangle^\ast$
  \item[(iv)]\ $\langle x, x\rangle\geq 0$
  \end{enumerate}
 \end{defn}

  The set of all adjointable maps $T:E\lra E$ will be denoted by $\cL(E)$. Since $\langle\cdot, \cdot\rangle$ is not strictly positive, the adjoint of a map from $\cL(E)$ is in general not unique. If $E$ is a \textit{Hilbert $A$-module}, that is the inequality at (iv) is strict and $E$ is complete with respect to the norm $\xi\mapsto ||\langle \xi,\xi\rangle ||^{\frac{1}{2}}$, then the adjoint is unique and the bounded elements of $\cL(E)$ form a C$^\ast$-algebra.

  If $\cB\subseteq\cA$, $\cB\subseteq \cD$ are unital inclusions of C$^\ast$-algebras, $\phi:\cA\lra\cD$ is a unital positive $\cB$-bimodule map and $a$ is a selfadjoint element of $\cA$, we will denote by $\cB\langle a \rangle$ the $\ast$-algebra generated in $\cA$ by $\cB$ and $a$ and by $\mu_a$, ``the $\cD$-distribution'' of $a$, that is the positive $\cB$-bimodule map $\phi_a:\bx\lra \cD$ defined by
 $\phi_a=\phi\circ \tau_a$
 where $\tau_a:\bx\lra \cA$ is the unique homomorphism such that $\tau_a(\X)=a$ and $\tau_a(b)=b$ for all $b\in\cB$.
  The set of elements from $\Sbd$ that can be realized in such way is $\Sbd^0$, more precisely we have the following property:

 \begin{prop}\label{prop1.2}
   Let $\mu\in\Sbd$. Then $\mu\in\Sbd^0$ if and only if there exist a C$^\ast$-algebra $A$ containing $\cB$ as a C$^\ast$-subalgebra, a completely positive $\cB$-bimodule map $\phi:A\lra \cD$ and a self-adjoint element $a\in A$ such that $\mu=\phi_a$.

   Moreover, the condition \emph{(\ref{pr2})} is equivalent to the existence of $M>0$ such that for  all $b_1, \dots, b_n\in M_m(\cB)$ we have
 \begin{equation}\label{pr2.b} ||(\id_{m}\otimes\mu)(\X\cdot b_1 \X\cdot b_2\cdots \X\cdot b_n \X)||< M^{n+1}||b_1||\cdots ||b_n||.
 \end{equation}
 where $\X$ acts on $M_m(\mathbb{C})\otimes\bx$ by multiplication on each entry (that is we identify $\X$ to $\id_m\otimes\X$).
 \end{prop}

\begin{proof}
If $\mu=\phi_a$ as above, then the result is trivial. 

Suppose now that $\mu\in\Sbd^0$. We first prove that $\mathcal{N}_0=\{f\in\bx : \ \mu(f^\ast f)=0\}$ is a left ideal of $\bx$. It suffices to prove that if $f\in\mathcal{N}_0$ then $\X\cdot f\in\mathcal{N}_0$ and $b\cdot f\in\mathcal{N}_0$ for all $b\in\cB$.
 
Since $b^\ast b\leq \|b^\ast b\|$ (in the C$^\ast$-algebra $\cB$), the positivity of $\mu$ implies
\[ \mu(f^\ast (\|b^\ast b\|-b^\ast b)f)\geq 0,\]
that is $\mu(f^\ast b^\ast b f)\leq \|b^\ast b\| \mu(f^\ast f)=0$, hence $b\cdot f\in\mathcal{N}_0$.

For $g=b_0 \X b_1\X\cdots \X b_n$ a monomial in $\bx$, define $\mathfrak{p}(g)=M^n\|b_0\|\|b_1\|\cdots\|b_n\|$  (in particular, condition  (\ref{pr2}) states that $\|\mu(g)\|\leq\mathfrak{p}(g)$). Consider 
\[
\cB\langle\langle \X\rangle\rangle_\mu=\{\sum_{n=0}^\infty f_n : f_n=\text{monomials in}\ \bx \text{such that}\ \sum_{n=0}^\infty \mathfrak{p}(f_n)<\infty\}.
\]
$\cB\langle \langle \X\rangle\rangle_\mu$ is a $\ast$-algebra (with the structure inherited from $\bx$) and $\mu$ extends to a positive map $\widetilde{\mu}:\cB\langle \langle \X\rangle\rangle_\mu\lra \cD$ via  $\widetilde{\mu}(\sum_{n=0}^\infty f_n)=\sum_{n=0}^\infty \mu(f_n)$.

Take now $g_n=(2n)!\left[(1-2n)(n!)^2 (4M)^{2n}\right]^{-1}\X^{2n}$, $n\geq 0$. Then $g_n=g_n^\ast$ and $\mathfrak{p}(g_n)\leq 4^{-n}$, so $g=\sum_{n=0}^\infty g_n$ is also a selfadjoint element of $\cB\langle \langle \X\rangle\rangle_\mu$. 

Since $g^2=1-\left[(2M)^{-1}\X\right]^2$, we have that 
\begin{align*}
0\leq\widetilde{\mu}(f^\ast g^\ast g f)&=\widetilde{\mu}\left(f^\ast [1-(2M)^{-2}\X^2]\cdot f\right)\\
&=\mu(f^\ast f)-(2M)^{-2}\mu\left((\X f)^\ast \X f\right)
\end{align*}
hence $\mu\left((\X f)^\ast \X f\right)\leq 4M^2 \mu(f^\ast f)$, so $\X f \in\mathcal{N}_0$.

 Consider now $\cK=\bx\otimes_\cB\cD$ with the right $\cD$-module structure given by
 $(f(\X)\otimes d)d_1=f(\X)\otimes dd_1$; also, since $\mu$ satisfies the condition (\ref{prop}), we can define a $\cD$-valued sesquilinear inner-product structure on $\cK$ (see \cite{lance}, page 40) via
 \[
 \langle f(\X)\otimes d_1, g(\X)\otimes d_2\rangle= d_2^\ast\mu(g(\X)^\ast f(\X) )d_1.
 \]
Let now $\mathcal{N}=\{ \eta\in\cK \colon \langle \eta, \eta\rangle =0\}$. From  the above argument on $\mathcal{N}_0$, we have that $\mathcal{N}$  is a left $\bx$-module.
Finally, take $E$ the completion of $\cK\slash \mathcal{N}$ in the norm induced by the inner-product structure (see \cite{lance}) and let $\xi=1\otimes 1+\mathcal{N}$. The multipliers with polynomials from $\bx$ are in the C$^\ast$-algebra of bounded maps from $\cL(E)$, since  condition (\ref{pr2}) ensures the boundness. Moreover, if $\phi(\cdot)=\langle \cdot \xi, \xi\rangle$ and $a$ is the right multiplier with $\X$, then $\mu=\phi_a$.

Condition (\ref{pr2.b}) is implied by the equality $\| \id_m\otimes a\| =\|a\|$, where the first norm is in the C$^\ast$-algebra
$M_m(A)$ and the second norm in the C$^\ast$-algebra $A$.

\end{proof}

%%%%%%%%%%%%%%%%%%%%%%%%%%%%%%%%%%%%%%%%%%%%%%%%%%%%%%%%%%%%
%%%%%%%%%%%%%%%%%%%%%%%%%%%%%%%%%%%%%%%%%%%%%%%%%%%%%%%%%%%%%
%%%%%%%%%%%%%%%%%%%%%%%%%%%%%%%%%%%%%%%%%%%%%%%%%%%%%%%%%%%%%%%%%%%%%%%%%%%%%%%%%%%%

 \section{Infinite divisibility: the boolean case}

The main result of this section is the non-commutative analogue of the following theorem (see \cite{speicherwaroudi}):
 \begin{thm}Any compactly supported real measure is infinitely divisible with respect to boolean convolution.
 \end{thm}
 We will use the following notion of boolean independence over a C$^\ast$-algebra $\cB$ (see \cite{mvcollmath}):
 \begin{defn}\label{defn1}
  Let $B$ be a unital C$^\ast$-algebra, $B\subseteq D$, $B\subseteq A$ be unital inclusions of $\ast$-algebras and $\phi: A\lra D$ be a unital $B$-bimodule map. A family $\{a_i\}_{i\in I}$
  of selfadjoint elements from $ A$ is said to be \emph{boolean independent} with respect to  $\phi$ if
  \[
  \phi(A_1 A_2 A_3\cdots)=\phi(A_1) \phi(A_2) \phi (A_3)\cdots
  \]
   for all $A_k\in\cB\langle a_{\epsilon(k)} \rangle_0$ (the $\ast$-algebra spanned by non-commutative polynomials in $a_{\epsilon(k)}$ and coefficients in $\cB$ without a free term), with $\epsilon(k)\in I$ and $\epsilon(k)\neq\epsilon(k+1)$.
 \end{defn}

  Let now $N\in\mathbb{N}$ and $\{\mu_j\}_{j=1}^N$ be  a family of elements from $\Sbd$. We define their  additive boolean convolution as follows. Consider the symbols $\{X_j\}_{j=1}^N$  such that $\mu_{X_j}:\cB\langle X_j\rangle\cong \bx\lra\cB$ coincides to $\mu_j$. Then consider the $\ast$-algebra $\cB\langle X_1, X_2, \dots, X_N \rangle$ with the conditional expectation $\mu$ such that its restrictions to $\cB\langle X_j\rangle$ coincide to $\mu_j$ and the mixed moments of $X_1, X_2, \dots, X_N$ are calculated via the rule from Definition \ref{defn1}. The additive boolean convolution of $\{\mu_j\}_{j=1}^N$ is the unital  $\cB$-bimodule map
  \[\uplus_{j=1}^N\mu_j=\mu_{X_1+X_2+\cdots+X_N}:\cB\langle X_1+X_2+\dots+X_N\rangle\cong\bx\lra\cD\]

  \begin{remark}\label{rem01} If $\mu_j\in\Sbd$ for al $j=1, \dots, N$, then $\uplus_{j=1}^N\mu_j$ is also in $\Sbd$.
  \end{remark}
  \begin{proof}
   For each $j$, consider the right-$\cD$ module $\cK_j=\cB\langle X_j\rangle\otimes_\cB\cD$ as in Proposition \ref{prop1.2}.

 Note that the $\cD$-submodule $1\otimes\cD$ is complemented in each $\cK_j$, since, for all $f(X)\in\bx$, we have that $1\otimes \mu_j(f(X_j))\in1\otimes_\cD$ and
 \[
 \langle 1\otimes1, f(X_j)\otimes 1-1\otimes \mu_j(f(X_j))\rangle_j =\mu_j(f(X_j)^\ast)-\mu_j(f(X_j))^\ast=0
 \]

 Denote by $\cK_j^0$ the complement of $1\otimes_\cD$ in $\cK_j$ and let
 \[ \cK=(1\otimes \cD)\oplus \bigoplus_{j=1}^N \cK_j^0.\]
 On $\cL(\cK)$ we consider the map $\phi:T\mapsto \langle T 1\otimes 1, 1\otimes 1\rangle$.

  $\cB\langle X_j\rangle$ can be seen as a algebra of linear maps on $\cK_j$ via
  \[
  f(X_j)\bigl[g(X_j)\otimes d\bigr]=\bigl(f(X_j)g(X_j)\bigr)\otimes d. \]
  Since $f(X_j)^\ast$ is adjoint to $f(X_j)$ we have $\cB\langle X_j\rangle\in\cL(\cK)$. moreover, by setting the restrictions of $\cB\langle X_j\rangle _{0}$ to each $\cK_l^0$ to the 0 if $l\neq j$, we can see $\cB\langle X_j\rangle$ as a subalgebra of $\cL(\cK)$. Note that $\phi_{|\cB\langle X_j\rangle_0}=\mu_j$.

  If $j\neq l$ and $A_1\in\cB\langle X_j\rangle_0$ while $A_2\in \cB\langle X_l\rangle_0$, then
  \begin{eqnarray*}
  A_1A_2 (1\otimes1)&=&A_1[\phi(A_2)(1\otimes1)+\xi], \ \text{where}\ \xi=A_2 (1\otimes1)-\phi(A_2)1\otimes1\in\cK_l^0\\
  &=&A_1\phi(A_2)(1\otimes1).
  \end{eqnarray*}
 Iterating, we obtain that for $A_k\in\cB\langle X_{\epsilon(k)} \rangle_0$  with $\epsilon(k)\neq\epsilon(k+1)$
 \[
 \phi(A_1\cdots A_m)=\phi(A_1\phi(A_2)\cdots \phi(A_m))=\phi(A_1)\cdots \phi(A_m).
 \]
 that is $X_j$'s are boolean independent in $\cL(\cK)$ with respect to $\phi$, and since the restrictions of $\phi$ to $\cB\langle X_j\rangle$ are $\mu_j$, we have that
 $
 \uplus_{j=1}^N\mu_j=\phi_{|\cB\langle X_1+\dots+X_N\rangle}$
 so q.e.d..

  \end{proof}

 \begin{defn}\label{defnbi1}
  An element $\mu\in\Sbd$ is said to be $\biguplus$-infinite divisible if for any positive integer $n$ there exist some $\mu_n\in\Sbd$ such that $\mu$ is the additive boolean convolutions of $n$ copies of $\mu_n$.
 \end{defn}
  The main result of this section si the following
 \begin{thm}\label{main1}Any element $\mu$ from $\Sbd$ is boolean infinitely divisible.
\end{thm}

  To prove \ref{main1}, we will need the following generalization of the notion of scalar boolean cumulants from \cite{speicherwaroudi}:
 \begin{defn}\label{defn23}
 In the above setting, let $X$ be a selfadjoint element from $\cA$. The \emph{boolean cumulants} of $X$ are defined as the multilinear maps $\{B_{n, X}\}_{n\geq1}$, with $B_{n, X}:\cB^n\lra\cD$, given by the following recurrence:
 \[
 \phi(Xb_nXb_{n-1}\cdots Xb_1)=\sum_{k=1}^n\phi(Xb_{n}\cdots Xb_{k+1})B_{k, X}(b_k, \cdots, b_1).
  \]
 \end{defn}

 As shown in \cite{mvcollmath}, the $\cD$-valued boolean cumulants defined above have the same additivity property as their scalar analogues (see \cite{mvcollmath}, Corollary 4.6):

 \begin{prop}\label{boolcum}
  If, in the above setting, $X$ and $Y$ are boolean independent over $\phi$, then, for all positive integers $n$, we have that
  \[B_{n, X+Y}=B_{n, X}+B_{n, Y}.\]
 \end{prop}

  Definition \ref{defnbi1} can be reformulated in terms of Proposition \ref{boolcum}. Namely, for $\mu\in\Sbd$ we define the $n$-th boolean cumulant of $\mu$ as the multilinear map
$B_{n, \mu}:\cB^n\lra\cD$ given by the recurrence:
\[\mu(\X b_1 \X b_2\cdots \X b_n)=\sum_{k=1}^nB_{k, \mu}(b_1, \dots, b_k)\mu(\X b_{k+1}\cdots \X b_n).
\]
From Proposition \ref{boolcum},  we have the following
\begin{remark}\label{remark1}
 An element $\mu\in\Sbd$ is $\biguplus$-infinite divisible if for any positive integer $n$ there exist $\mu_n\in\Sbd$ such that for all positive integers $m$ we have that
 \[B_{m, \mu}=nB_{m, \mu_n}.\]
\end{remark}

Before proving the main result, i. e. Theorem \ref{main1}, we will first need the following result.

\begin{lemma}\label{lema1}
Let $B\subseteq A$ be a unital inclusion of C$^\ast$-algebras, $\cH$ be a semi-inner product $A$-bimodule which is also a left $B$-module and $\Omega$ a symbol that commutes with $A$ such that $\langle \Omega, \Omega\rangle=1$ and $\langle \cdot \Omega, \Omega\rangle:\mathcal{L}(\cH\oplus\Omega A)\lra A$ is a $B$-bimodule map.

Let $T, \Lambda\in \mathcal{L}(\cH\oplus \Omega A)$ be selfdajoint operators such that $T(\Omega A)=0, T(\cH)\subseteq \cH$, $\Lambda(\cH)=0$, $\Lambda(\Omega A)\subseteq \Omega A$. For  $\xi\in\cH$ define the operators $a_\xi, a^\ast_\xi\in\mathcal{L}(\cH+\oplus \Omega A)$ given by
\[\left\{
\begin{array}{cl}
a_\xi\alpha=0& \alpha\in A\\
a_{\xi}\eta=\langle\eta, \xi\rangle & \eta\in\cH
\end{array}
\right.
\hspace{2cm}
\left\{
\begin{array}{cl}
a^\ast_\xi\alpha=\xi\alpha& \alpha\in A\\
a^\ast_{\xi}\eta=0 & \eta\in\cH
\end{array}
\right..
\]
 Then $a_\xi$, $a^\ast_\xi$ are adjoint to each other and the boolean cumulants $\{B_{n, V}\}_n$ of $V=a_\xi+a^\ast_\xi+T+\Lambda$ with respect to the map $\langle\cdot \Omega,  \Omega\rangle$ are given by:
 \begin{eqnarray*}
 \text{(i)} &&B_{1, V}(b_1)=\langle \Lambda b_1\Omega, \Omega\rangle\\
 \text{(ii)}&& B_{n, V}(b_n, \dots, b_1)=
 \langle a_\xi b_n Tb_{n-1}\cdots Tb_2a^\ast_\xi b_1\Omega,
  \Omega\rangle\ \text{if}\ n\geq2.
 \end{eqnarray*}
\end{lemma}

\begin{proof}
 The fact that $a_\xi$ and $a^\ast_\xi$ are adjoint to each other is just a trivial computation. For (i), note that $b_1\Omega\in\Omega A$, hence $\langle Vb_1\Omega, \Omega\rangle=\langle\Lambda b_1\Omega, \Omega\rangle$.

 For (ii), remark first that
 \[a_\xi b_n Tb_{n-1}\cdots Tb_2a^\ast_\xi b_1\Omega\subset\Omega A\]
  for all $b_1, \dots, b_n\in B$, since $a^\ast_\xi b_1\Omega=\xi b_1\in\cH$, also $b_nTb_{n-1}\cdots Tb_2\cH\subseteq\cH$
  and $a_\xi\cH\subseteq\Omega A$.

  We have that
  \begin{eqnarray*}
  \langle Vb_n\cdots Vb_1\Omega, \Omega\rangle
  &=&
  \sum_{V_j\in\{a_\xi, a^\ast_\xi, T, \Lambda\}}
  \langle V_nb_n\cdots V_1b_1\Omega, \Omega\rangle.
  \end{eqnarray*}

  Let us suppose that the term $\langle V_nb_n\cdots V_1b_1\Omega, \Omega\rangle$ does not cancel. Since $a_\xi(\Omega A)=T(\Omega A)=0$, it follows that $V_1\in\{a^\ast_\xi, \Lambda\}$. If $V_1=a^\ast_\xi$, then $V_1b_1\Omega\in\cH$ and, since
  \[a^\ast_\xi\cH=\Lambda\cH=0\]
  it follows that $V_2\in\{a_\xi, T\}$.

  Also, since $T\cH\subseteq\cH$, if $V_2=V_3=\dots=V_p=T$, then $V_{p+1}\in\{a_\xi, T\}$. Finally, note that
 $ Tb_n\cdots Tb_2a^\ast\xi b_1\Omega\in\cH$,
  henceforth
  \begin{eqnarray*}
  \langle Vb_n\cdots Vb_1\Omega, \Omega\rangle
  &=&
  \langle Vb_n\cdots Vb_2\Lambda b_1\Omega, \Omega\rangle \\
  &&\hspace{1cm}+ \sum_{p=2}^n\langle Vb_n\cdots Vb_{p+1}
  a_\xi b_p Tb_{p-1}\cdots Tb_2a^\ast\xi b_1 \Omega, \Omega\rangle\\
  &=&\langle Vb_n\cdots Vb_2\Omega, \Omega\rangle\cdot\langle\Lambda b_1\Omega, \Omega\rangle\\
  &&\hspace{-.2cm}+ \sum_{p=2}^n\langle Vb_n\cdots Vb_{p+1}\Omega, \Omega\rangle\cdot\langle
  a_\xi b_p Tb_{p-1}\cdots Tb_2a^\ast\xi b_1 \Omega, \Omega\rangle\\
  \end{eqnarray*}

and the conclusion follows from (i) and Definition \ref{defn23}.
\end{proof}

\textit{Proof of the Theorem \ref{main1}}

 Consider $\cK=\bx\otimes_\cB\cD$ as for Remark \ref{rem01}. Denote by $\cK^0$ the complement of $1\otimes_\cD$ in $\cK$ and let $\xi=\X\otimes 1-1\otimes \mu(\X) \in \cK^0$. Define the operator $a_\xi:\cK\lra\cK$ by
 \begin{eqnarray*}
 a_\xi(1\otimes1)&=&0\\
 a_\xi\eta&=& 1\otimes\langle \eta, \xi\rangle\ \text{for}\ \eta\in\cK^0
 \end{eqnarray*}
 Remark that $a_\xi$ is $\cD$-linear, adjointable, and its adjoint is given by
  \begin{eqnarray*}
 a^\ast_\xi(1\otimes 1)&=&\xi 1\\
 a^\ast_\xi\eta&=&0\ \text{for}\ \eta\in\cK^0
 \end{eqnarray*}

Let also let $\textgoth{L}_\X:\cK\lra\cK$ be the selfadjoint map given by $ f(\X)\otimes d\mapsto \X f(\X)\otimes d$. We will identify $\cB$ with a subalgebra of $\mathcal{L}(\cK)$ via
\[b[f(\X)\otimes d]=[bf(\X)]\otimes d.\]
Note that $\langle \cdot 1\otimes1, 1\otimes1\rangle$ is a $\cB$-bimodule map.

We have that
\begin{eqnarray*}
\langle \textgoth{L}_\X b_1\cdots \textgoth{L}_\X b_n(1\otimes 1),  1\otimes 1\rangle
&=&
\langle  \X b_1\X b_2\cdots \X b_n\otimes, 1\otimes1\rangle\\
&=&\mu(\X b_1\X b_2\cdots \X b_n)
\end{eqnarray*}
hence the distribution of $\textgoth{L}_X$ with respect to $\langle \cdot 1\otimes1, 1\otimes1\rangle$ coincides with $\mu$.

Consider the selfadjoint map $\Lambda_\mu\in\mathcal{L}(\cK)$ given by
\[
\left\{
\begin{array}{ll}
\Lambda_\mu(1\otimes d)=1\otimes\mu(\X)d,& d\in\cD\\
\Lambda_\mu\eta=0, & \eta\in\cK^0
\end{array}
\right.
\]
and define $T\in\mathcal{L}(\cK)$ as the map $T=\textgoth{L}_\X-(a_\xi+a^\ast_\xi)-\Lambda_\mu$. We have that $T$ is selfadjoint, $T\cK^0\subseteq\cK^0$ and $T(1\otimes\cD)=0$.

 Since $\textgoth{L}_\X=a_\xi+a^\ast_\xi +\Lambda_\mu+T$ and its distribution with respect to $\langle \cdot 1\otimes1, 1\otimes1\rangle$ is $\mu$  from Lemma \ref{lema1} we have that the boolean cumulants of $\mu$ are given by
 \begin{eqnarray}
  &&B_{1, \mu}(b_1)=\langle \Lambda_\mu b_1(1\otimes1), 1\otimes1\rangle \label{bool1}\\
 && B_{n, \mu}(b_n, \dots, b_1)=\langle a_\xi b_n Tb_{n-1}\cdots Tb_2a^\ast_\xi b_1(1\otimes1), 1\otimes1\rangle\ \text{if}\ n\geq2. \label{bool2}
 \end{eqnarray}

  Fix $N$ a positive integer. Let $\xi_N=\frac{1}{\sqrt{N}}\xi$, $\Lambda_N=\frac{1}{N}\Lambda_\mu$ and $Y_N\in\mathcal{L}(\cK)$ be the selfadjoint operator
  \[Y_N=a_{\xi_N}+a^\ast_{\xi_N}+\Lambda_N+T.\]

  Define $\mu_N\in\Sbd$ via
  \[\mu_N(f(X))=\langle f(Y_N)(1\otimes1), 1\otimes 1\rangle.\]

  From Lemma \ref{lema1}, the boolean cumulants of $\mu_N$ are given by:
  \begin{eqnarray*}
  B_{1, \mu_N}(b_1)
  &=&
  \langle \Lambda_N b_1(1\otimes1), 1\otimes1\rangle
  =
  \frac{1}{N}\langle \Lambda_\mu b_1(1\otimes1), 1\otimes1\rangle\\
  &=&
  \frac{1}{N} B_{1, \mu}(b_1), 
  \end{eqnarray*}

  \begin{eqnarray*}
  B_{n, \mu_N}(b_1, \dots, b_n)
  &=&
  \langle a_{\xi_N} b_1 Tb_{2}\cdots Tb_{n-1}a^\ast_{\xi_N} b_n(1\otimes1), 1\otimes1\rangle\\
  &=&
  \frac{1}{N}\langle a_{\xi} b_1 Tb_{2}\cdots Tb_{n-1}a^\ast_{\xi} b_n(1\otimes1), 1\otimes1\rangle\\
  &=&\frac{1}{N}B_{n, \mu}(b_1, \dots, b_n)
  \end{eqnarray*}
  and the conclusion follows from Remark \ref{remark1}.

  \hfill $\square$

%%%%%%%%%%%%%%%%%%%%%%%%%%%%%%%%%%%%%%%%%%%%%%%%%%%%%%%%%%%%%%%%%%%%%%%%%
%%%%%%%%%%%%%%%%%%%%%%%%%%%%%%%%%%%%%%%%%%%%%%%%%%%%%%%%%%%%%%%%%%%%%%%%%%%%
%%%%%%%%%%%%%%%%%%%%%%the free case%%%%%%%%%%%%%%%%%%%%%%
%%%%%%%%%%%%%%%%%%%%%%%%%%%%%%%%%%%%%%%%%%%%%%%%%%%%%%%%%%%%%%%%%%%%%%%%%%

\section{Infinite divisibility: the Free case}\label{section3}

%%%%%%%%%%%%%%%%%%%%%%%%%%%%%%%%%%%%%%%%%%%%%%%%%%%%%%%%%%%%%%%%%%%%%
%%%%%%%%%%%%%%%%%%%%%%%%%%%%%%%%%%%%%%%%%%%%%%%%%%%%%%%%%%%%%%%%%%%%%%%%%%%%
%%%%%%%%%%%%%%%%%%%%%%%%%%%%%%%%%%%%%%%%%%%%%%%%%%%%%%%%%%%%%%%%%

\subsection{Preliminaries}\label{section31}

 \begin{defn}\label{defn31}\emph{(see \cite{voi95})}
  Let $\cB$ be a unital C$^\ast$-algebra,  $\cB\subseteq\cA$ be a  unital inclusion of $\ast$-algebras and $\phi:\cA\lra\cB$ be a positive conditional expectation. A family $\{X_i\}_{i\in I}$
  of selfadjoint elements from $\cA$ is said to be \emph{free} with respect to if $\phi$ if
 \[\phi(A_1A_2\cdots A_n)=0\]
  whenever $A_j\in\cB\langle X_{\epsilon(j)}\rangle$ with $\epsilon(k)\neq \epsilon(k+1)$ and $\phi(A_j)=0$.
 \end{defn}

  In the above setting, let $X$ be a selfadjoint element from $\cA$. The free cumulants of $X$ are the multilinear functions $\kappa_{n, X}:\cB^n\lra\cB$ given by the recurrence:
  \begin{eqnarray*}
  \phi(Xb_1Xb_2\cdots Xb_n)&=& \\
  &&\hspace{-2.3cm}\sum_{p=1}^n\sum_{1<j_1<\dots<j_{p-1}}\kappa_{p, X}((b_1\phi(Xb_2\cdots X b_{j-1}), (b_{j_1+1}\phi(Xb_{j_1+2}\cdots Xb_{j_2}), \dots \\
  &&\hspace{2cm}\dots, (b_{j_{p_1}}\phi(Xb_{j_{p-1}+1}\cdots Xb_n)).
  \end{eqnarray*}

  More intuitively, the above relation can be graphically illustrated by the picture below, where the boxes stand for the application of $\phi$ and the oblique lines signify that each $Y_s=\phi(Xb_{j_s+2}\cdots Xb_{j_{s+1}})$ are multiplied with $b_{j_{s}+1}$ in the arguments of the free cumulants:

  \setlength{\unitlength}{.15cm}
 \begin{equation*}
 \begin{picture}(17,7)

 \put(-27,2){\line(0,1){3}}
 \put(-27,2){\line(1,0){12}}
\put(-15,2){\line(0,1){3}}

\put(-27,5){\line(1,0){12}}
\put(-14,3){=}

\put(-12,3){\huge{$\sum$}}

\put(-7,2){\line(0,1){4}}

\put(-7,6){\line(1,0){14}}

\put(-8,5){\line(3,-2){5}}

\put(0,2){\line(0,1){4}}

\put(-1,5){\line(3,-2){5}}

\put(14,5){\line(3,-2){5}}

\put(21,4){\line(3,-1){6}}

\put(22,2){\line(0,1){4}}

\put(14,6){\line(1,0){8}}

\put(15,2){\line(0,1){4}}

\qbezier[7](7,4)(9,4)(13,4)

\put(-6,2){\line(0,1){3}} \put(-2,2){\line(0,1){3}}

\put(-6,2){\line(1,0){4}}

\put(-6,5){\line(1,0){4}}

\put(1,2){\line(0,1){3}} \put(5,2){\line(0,1){3}}

\put(1,2){\line(1,0){4}}

\put(1,5){\line(1,0){4}}

\put(16,2){\line(0,1){3}} \put(20,2){\line(0,1){3}}

\put(16,2){\line(1,0){4}}

\put(16,5){\line(1,0){4}}

\put(23,2){\line(0,1){3}}

\put(23,2){\line(1,0){9}}

\put(32,2){\line(0,1){3}}

\put(23,5){\line(1,0){9}}

 \end{picture}
 \end{equation*}

  The free cumulants have the following additivity property (see \cite{speicherhab}, \cite{mvcomstoc}):
  \begin{prop}\label{prop32}If $X, Y$ are free in the sense of Definition \ref{defn31}, then
  \[\kappa_{n, X+Y}=\kappa_{n, X}+\kappa_{n, Y}.\]
  \end{prop}

  Let now $N\in\mathbb{N}$ and $\{\mu_j\}_{j=1}^N$ be  a family of elements from $\Sigma_\cB$. We define their additive free  convolution similarly to be boolean case:  Consider the symbols $\{X_j\}_{j=1}^N$; on the algebra $\cB\langle X_1, X_2, \dots, X_N \rangle$ take the conditional expectation $\mu$ such that $\mu\circ \tau_{X_j}=\mu_j$ and the mixed moments of $X_1, \dots, X_n$ are computed via the rules  from Definition \ref{defn31}. The free additive convolution of $\{\mu_j\}_{j=1}^N$ is the conditional expectation
  \[\boxplus_{j=1}^N \mu_j=\mu\circ\tau_{X_1+X_2+\cdots+X_N}:\cB\langle X_1+X_2+\dots+X_N\rangle\cong\bx\lra\cB.\]

 We have that $\boxplus_{j=1}^N \mu_j$ is also an element of $\Sb$: in \cite{speicherhab}, Theorem 3.5.6, it is shown that $\mu$, defined as above, is a positive conditional expectation, therefore so is $\mu\circ\tau_{X_1+X_2+\cdots+X_N}$.

 \begin{defn}
  An element $\mu\in\Sigma_\cB$ is said to be $\boxplus$-infinite divisible if for any positive integer $N$ there exists $\mu_N\in\Sigma_\cB$ such that $\mu$ is the free additive  convolution of $N$ copies of $\mu_N$.
 \end{defn}

 \subsection{Free cumulants and $\boxplus$-infinite divisibility}

 \begin{defn}\label{def33}
  Let $\mu\in\Sb$. Using the relations from Section \ref{section31}, we define the free cumulants of $\mu$ as the multilinear maps $\kappa_{n, \mu}:\cB^n\lra\cB$ given by the recurrence:
   \begin{eqnarray*}
  \mu(\X b_1 \X b_2\cdots \X b_n)&=&\\
  &&\hspace{-2cm}\sum_{p=1}^n\sum_{1<j_1<\dots<j_{p-1}}\kappa_{p, \mu}((b_1\mu(\X b_2\cdots \X b_{j-1}), (b_{j_1+1}\mu(\X b_{j_1+2}\cdots \X b_{j_2}), \dots\\
  &&\hspace{2cm}\dots, (b_{j_{p_1}}\mu(\X b_{j_{p-1}+1}\cdots \X b_n)).
  \end{eqnarray*}
 \end{defn}
 \begin{remark}\label{remark33}
 Using Proposition \ref{prop32}, we can reformulate Definition \ref{def33} in terms of free cumulants. More precisely, $\mu\in\Sb$ is $\boxplus$-infinitely divisible if for any positive integer $N$ there exists some $\mu_N\in\Sb$ such that for all $m$
 \[\kappa_{m,\mu}=N\cdot\kappa_{m, \mu_N}.\]
 \end{remark}

 For $\mu\in\Sb$, define the conditional expectation $\rho_\mu:\bx\lra\cB$ generated by
 \[\rho_\mu(\X b_1\X b_2\cdots \X b_n)=\kappa_{n, \mu}(b_1, b_2, \cdots, b_n).\]

 \begin{prop}\label{rhoprop}
  If $\mu\in\Sb$ is $\boxplus$-infinitely divisible then the restriction of $\rho_\mu$ to $\bx_0$ is positive.
 \end{prop}
 \begin{proof}
  Fix $N$ and suppose that $\mu$ is the free additive convolution of $N$ copies of $\mu_N$.
  Note that, for $n\leq 1$,
  \begin{equation}\label{o1}
  \mu_N(\X b_1\cdots {\X}b_n)=\frac{1}{N}\kappa_{n, \mu}(b_1, \dots, b_n)+O(\frac{1}{N^2}).
  \end{equation}

  The assertion is trivial for $n=1$. Suppose that (\ref{o1}) is true for $n<m$. Since the free cumulants are multilinear, for all $1=l_1<l_2<\dots<l_{p+1}=m+1$ and $Y_s=b_{l_s}\mu_N({\X}b_{l_s+1}\cdots {\X}b_{l_{s+1}-1})$, ($1\leq s \leq p$) we have that
  \[
    Y_s=\left \{
    \begin{array}{ll}
     b_{l_s}& \text{if}\ l_{s+1}=l_s+1\\
    O(\frac{1}{N})& \text{if}\ l_{s+1}=l_s+1
    \end{array}
    \right.
  \]
  hence

  $\displaystyle
   \kappa_{p,\mu_N}(Y_1, \dots, Y_p)
   =\frac{1}{N}\kappa_{p, \mu}(Y_1, \dots, Y_p)
   =O(\frac{1}{N^2})
   $\ unless $l_{s+1}=l_s+1$ for all $s\in\{1, \dots, p\}$, i. e. $p=m$.

   Definition \ref{def33} gives
   \[
    \mu_N({\X}b_1\cdots {\X}b_m)=\kappa_{m, \mu_N}(b_1, \dots, b_m)+O(\frac{1}{N^2}) = \frac{1}{N}\kappa_{m, \mu}(b_1, \dots, b_m)O(\frac{1}{N^2})\]
    that is (\ref{o1}). It follows that
    \begin{equation}\label{eq32}
    \lim_{N\lra\infty} N\cdot\mu_N({\X}b_1\cdots {\X}b_n)=\kappa_{n, \mu}(b_1, \dots, b_n)
    \end{equation}
  Fix now $f(\X)\in\bx_0$. Then
    \[
    {}\rho_{\mu}(f({\X})^\ast f({\X}))
    =\lim_{N\lra\infty}[N\mu_N(f({\X})^\ast f({\X}))]\geq 0\ \text{since all $\mu_N$ are positive}.
    \]
 \end{proof}

 \begin{lemma}\label{freecumullemma}
  Let $A$ be a unital C$^\ast$-algebra and $\cH$ a semi-inner-product $A$-module which is also a left $A$-module. Let $\cT({\cH})$ be the \emph{full Fock $A$-module} over $\cH$, that is
  \[\cT(\cH)=A\oplus \cH\oplus (\cH\otimes_{B}\cH)\oplus(\cH\otimes_{B}\cH\otimes_{B}\cH)\oplus\dots.\]

  Fix $\xi\in\cH$, $\beta\in A$, $T\in\mathcal{L}(\cH)$, $T$ and $\beta$ selfadjoint, and define the maps ($\alpha\in A$, $\eta_1,\dots,\eta_n\in\cH$)
  \begin{eqnarray*}
  &&\left\{
\begin{array}{l}
a_\xi\alpha=0\\
a_{\xi}\eta_1\otimes\eta_2\otimes\cdots\otimes\eta_n=\langle\eta,
\xi\rangle\eta_2\otimes\cdots\otimes\eta_n
\end{array}
\right.\\
&&
\left\{
\begin{array}{l}
a^\ast_\xi\alpha=\xi\alpha\\
a^\ast_{\xi}\eta_1\otimes\cdots\otimes\eta_n=\xi\otimes\eta_1\otimes\cdots\otimes\eta_n
\end{array}
\right.\\
&&
\left\{
\begin{array}{l}
p(T)\alpha=0\\
p(T)\eta_1\otimes\eta_2\otimes\cdots\otimes\eta_n=T(\eta_1)\otimes\eta_2\otimes\cdots\otimes\eta_n
\end{array}
\right.
\end{eqnarray*}

  Then $a_\xi$, $a^\ast_\xi$ are adjoint to each other, $p(T)$ is selfadjoint and the cumulants of $X=a_\xi+a^\ast_\xi+\widetilde{T}+\beta\chi_A$ with respect to the conditional expectation $\langle \cdot 1, 1\rangle$ are given by
  \begin{eqnarray*}
  \kappa_{1, X}(b_1)&=&\beta b_1\\
  \kappa_{n, X}(b_1, \dots, b_n)&=&\langle a_\xi b_1 p(T)b_2\cdots p(T)b_{n-1}a^\ast\xi b_1 1, 1\rangle.
  \end{eqnarray*}
 \end{lemma}

 \begin{proof}

   For $\kappa_{1, X}$ the assertion is trivial. To prove the relation for higher order free cumulants, let us fix $N>1$ and consider $\{(\cH_{j,N}, T_{j, N}, \xi_{j,N})\}_{j=1}^N$ to be a set of $N$ identical copies of $(\cH, T, (\frac{1}{\sqrt{N}})\xi)$ from above. Let
   \begin{eqnarray*}
   \cH_N
   &=&
   \cH_{1, N}\oplus\cH_{2, N}\oplus\dots\oplus\cH_{N,N}\\
   \xi_N
   &=&
   \xi_{1,N}\oplus\xi_{2, N}\oplus\dots\oplus\xi_{N, N}\in\cH_N\\
   T_N
   &=&
   T_{1, N}\oplus\dots \oplus T_{N, N}\in\cL(\cH_N)\\
   X_N&=&a_{\xi_N}+a_{\xi_N}^\ast+P(T_N)+\beta\chi_A\in\cL(\cH_N)\\
   X_{j, N}
   &=&a_{\xi_{j, N}}+a_{\xi_{j, N}^\ast}+p(T_{j, N})+\frac{1}{N}\beta\chi_A\in\cL(\cH_N).
   \end{eqnarray*}

 First, note that $X$ and $X_N$ are identically distributed with respect to $\langle \cdot 1, 1\rangle$. Then $X_N=X_{1, N}+\dots+ X_{N, N}$, $X_{j, N}$ are identically distributed and free. To see that, consider $A_1, \dots, A_m$, $A_k\in A\langle X_{\epsilon(k), N}\rangle$ with $\langle A_k1, 1\rangle=0$ and $\epsilon(k)\neq\epsilon(k+1)$. We have to show that $\langle A_m\cdots A_11, 1\rangle=0$.

 Since $\langle A_1 1, 1\rangle=0$, we have that $A_11\in\cT(\cH_{\epsilon(1), N})\ominus A$. Also, $\langle A_2 1, 1\rangle=0$, so $A_2A_1 1\in(\cT(\cH_{\epsilon(2), N})\ominus A)\otimes(\cT(\cH_{\epsilon(1), N})\ominus A)$. Iterating, we obtain
  \[
  A_m\cdots A_2A_1 1\in (\cT(\cH_{\epsilon(m), N}\ominus A)\otimes\cdots\otimes(\cT(\cH_{\epsilon(1), N})\ominus A)
  \]
so $\langle A_m\cdots A_2A_1 1, 1\rangle=0$.

Using (\ref{eq32}) and denoting $\mathcal{V}=\{\frac{1}{\sqrt{N}}a_\xi, \frac{1}{\sqrt{N}}a^\ast_\xi, p(T), \frac{1}{N}\beta\chi_A\}$ we have ($n\geq2$):
\begin{eqnarray*}
\kappa_{n, X}(b_1, \dots, b_n)
&=&
 \kappa_{n, X_{N}}(b_1, \dots, b_n)\\
&=&
\lim_{N\lra\infty}N\langle X_Nb_1\cdots X_Nb_n 1, 1\rangle\\
&=&\sum_{V_i\in\mathcal{V}}\langle V_1b_1\cdots V_nb_n 1, 1\rangle
\end{eqnarray*}

But $\langle V_1b_1\cdots V_nb_n 1, 1\rangle=0$ unless $V_n\in\{\frac{1}{\sqrt{N}}a^\ast_\xi, \frac{1}{N}\beta\chi_A\}$ and $V_1\in\{\frac{1}{\sqrt{N}}a_\xi, \frac{1}{N}\beta\chi_A\}$, therefore
\[
\langle X_Nb_1\cdots X_Nb_n 1, 1\rangle=\langle \frac{1}{\sqrt{N}}a_\xi b_1 p(T)b_2\cdots p(T)b_{n-1}\frac{1}{\sqrt{N}}a^\ast_\xi b_1 1, 1\rangle + O(\frac{1}{N^2})
\]
hence the conclusion.
 \end{proof}

 \begin{thm}\label{rhothm}
 The conditional expectation $\mu\in \Sb$ is $\boxplus$-infinitely divisible if and only if the restriction of  $\rho_{\mu}$ to $\bx_0$ is positive.
 \end{thm}

 \begin{proof}

  Suppose that $\rho_{\mu|\cB\langle X\rangle_0}$ is positive. Then (see \cite{lance}, pag. 42) $\bx_0$ is a semi-inner-product $\cB$-module with respect to the pairing
  \[\langle f(\X), g(\X)\rangle = \rho_\mu(g(\X)^\ast f(\X)).\]

  Consider the selfadjoint map $T:\bx_0\lra\bx_0$ given by $T(f(\X))=\X f(\X)$ and denote by $V$ the map from $\cL(\cT(\bx_0))$ defined as
  \[ V=a_\X+a^\ast_\X+\widetilde{T}+\mu(\X)\id.\]
 From Lemma \ref{freecumullemma} we have that the free cumulants of $V$ with respect to $\langle \cdot 1, 1\rangle$ are given by
 \begin{eqnarray*}
 \kappa_{1, V}(b_1)&=&\mu({\X})(b_1)=\kappa_{1, \mu}(b_1)\\
 \kappa_{n, V}(b_1, \dots, b_n)
 &=&
 \langle a_{\X} b_1 \widetilde{T} b_2\cdots \widetilde{T} b_{n-1}a^\ast_{\X} b_n 1, 1\rangle=
 \langle b_1{\X}b_2\cdots {\X}b_n, {\X}\rangle\\ &=& \kappa_\mu(b_1, \cdots, b_n)
 \end{eqnarray*}

 Fix a positive integer $N$. From Remark \ref{remark33} it suffices to find a selfadjoint element from $\cL(\cT(\bx_0))$ whose free cumulants are proportional to the free cumulants of $V$ by a factor of $\frac{1}{N}$. Define
 \[ V_N=\frac{1}{\sqrt{N}}a_{\X}+\frac{1}{\sqrt{N}}a^\ast_{\X}+\widetilde{T}+\frac{1}{N}\mu({\X})\id.\]
 Applying again Lemma \ref{freecumullemma}, the free cumulants of $V_N$ are
\begin{eqnarray*}
 \kappa_{1, V_N}(b_1)&=&\frac{1}{N}\mu({\X})(b_1)=\frac{1}{N}\kappa_{1, V}(b_1)\\
 \kappa_{n, V_N}(b_1, \dots, b_n)
 &=&
 \langle \frac{1}{\sqrt{N}}a_{\X} b_1 \widetilde{T} b_2\cdots \widetilde{T} b_{n-1}\frac{1}{\sqrt{N}}a^\ast_{\X} b_n 1, 1\rangle\\&=& \frac{1}{N}\kappa_{n, V}(b_1, \cdots, b_n)
 \end{eqnarray*}

  The converse implication is Proposition \ref{rhoprop}.

 \end{proof}

 %%%%%%%%%%%%%%%%%%%%%%%%%%%%%%%%%%%%%%%%%%%%%%%%%%%%%%%%
 %%%%%%%%%%%%%%%%%%%%%%%%%%%%%%%%%%%%%%%%%%%%%%%%%%%%%%
 %%%%%%%%%%%%%%%%%%%%%%%%%%%%%%%%%%%%%%%%%%%%%%%%%%%%%

 \section{Infinite divisibility: the c-free case}\label{cfreeinf}

 In this section we aim to extend the results from Section \ref{section3} to the case $\mu\in\Sbd$. First, we will need a suitable definition for the free additive convolution of elements from $\Sbd$; in this setting, if $\cB$ is simply replaced by $\cD$ in Definition \ref{defn31}, the resulting relation does not uniquely determine the joint moments of $X_1, \dots, X_n$. As shown in \cite{boca90}, \cite{dykemablanchard}, a more suitable approach is the conditional freeness (see also \cite{mvcomstoc}, \cite{bsk}).

 \begin{defn}\label{def41}
 Let $\cB$ be a unital C$^\ast$-algebra $\cB\subseteq\cA$, $\cB\subseteq \cD$ be unital inclusions of $\ast$-algebras, $\varphi:\cA\lra\cB$ be a conditional expectation and $\theta:\cA\lra\cD$ be a unital $\cB$-bimodule map.

 The family $\{X_i\}_{i\in I}$ of selfadjoint elements from $\cA$ is said to be c-free with respect to $(\theta, \varphi)$ if
 \begin{enumerate}
 \item[(i)] the family $\{X_i\}_{i\in I}$ is free with respect to $\varphi$
 \item[(ii)] $\theta(A_1A_2\cdots A_n)=\theta(A_1)\theta(A_2)\cdots \theta(A_n)$ for all $A_i\in\cB\langle X_{\epsilon(i)}\rangle$ such that $\varphi(A_i)=0$ and $\epsilon(k)\neq\epsilon(k+1)$.
 \end{enumerate}
 \end{defn}

 Let $X$ be a selfadjoint element from $\cA$. The c-free cumulants of $X$ are the multilinear functions $\cka_{n, X}:\cB^n\lra\cD$ given by the recurrence:
  \begin{eqnarray*}
  \theta({\X}b_1{\X}b_2\cdots {\X}b_n)&=&
  \sum_{p=1}^n\sum_{\substack{1<j_1<\dots<j_{p}\\l_ p=n-1}}
  \theta({\X}b_1\cdots {\X}b_{l_1})\cdot\\
  &&\hspace{-2.6cm}\cka_{p, {\X}}\left(b_{l_1+1}\varphi({\X}b_{l_1+2}\cdots {\X}b_{l_2}),\dots, b_{l_{p-1}+1}\varphi({\X}b_{l_{p-1}+2}\cdots {\X}b_{l_p}), b_n\right)\\
  \end{eqnarray*}

  As in the previous section, the above equation can be represented more intuitively by the picture below, where the dark boxes stand for the application of $\theta$, the light ones for the application of $\varphi$ and the oblique lines signify that each $Y_s=\phi({\X}b_{l_s+2}\cdots {\X}b_{l_{s+1}})$ are multiplied with $b_{l_{s}+1}$ in the arguments of the c-free cumulants, except for $b_n$.

 %%%%%%%%%%%%%%%%%%%%%%%%%%%%%%%%%%

 \setlength{\unitlength}{.15cm}
 \begin{equation*}
 \begin{picture}(10,6)

 \put(-27,2){\rule{12\unitlength}{3\unitlength}}
 \put(-14,3){=}

\put(-12,3){\huge{$\sum$}}

\put(-7,2){\rule{9\unitlength}{3\unitlength}}

\put(2.6,3){$\times$}

\put(5,2){\line(0,1){4}}

\put(5,6){\line(1,0){13}}

\put(4,5){\line(3,-2){5}}

\put(6,2){\line(0,1){3}}

\put(10,2){\line(0,1){3}}

\put(12,2){\line(0,1){4}}

\put(11,5){\line(3,-2){5}}

%\put(16,6){\line(1,0){8}}

\put(13,2){\line(0,1){3}}

\put(17,2){\line(0,1){3}}

\put(13,2){\line(1,0){4}}

\put(13,5){\line(1,0){4}}

\put(6,2){\line(1,0){4}}

\put(6,5){\line(1,0){4}}

\qbezier[7](19,4)(21,4)(25,4)

\put(27,2){\line(0,1){4}}

\put(26,6){\line(1,0){8}}

\put(34,2){\line(0,1){4}}

\put(28,2){\line(1,0){4}}

\put(28,5){\line(1,0){4}}

\put(32,2){\line (0,1){3}}

\put(28,2){\line(0,1){3}}

\put(26,5){\line(3,-2){5}}

\put(33,5.5){\line(1,-2){2}}

 \end{picture}
 \end{equation*} %
 %%%%%%%%%%%%%%%%%%%%%%%%%%%%%%%%%%%%%%%%%%%%%%%%%%%%%%%%%%%%%%%%%%

  The c-free cumulants have the following additivity property (see \cite{mvcomstoc}, \cite{mlotk}):
  \begin{prop}\label{prop42}If $X, Y$ are c-free with respect to $(\theta, \varphi)$ in the sense of Definition \ref{def41}, then
  \begin{eqnarray*}
  \kappa_{n, X+Y}&=&\kappa_{n, X}+\kappa_{n, Y}\\
  \cka_{n, X+Y}&=&\cka_{n, Y}+\cka_{n, Y}
  \end{eqnarray*}
  where $\kappa_{n, X}$ is the $n$-th free cumulant of $X$ with respect to the conditional expectation $\varphi$.
  \end{prop}

  Let now $N$ be a positive integer and $\{(\mu_i, \nu_i)\}_{i=1}^N$ be  a family from $\Sbd\times\Sb$. We define their additive c-free  convolution similarly to be boolean and free case:  Consider the selfadjoint symbols $\{X_i\}_{i=1}^N$  and the mappings
  \begin{eqnarray*}
  \mu:\cB\langle X_1, X_2, \dots, X_N \rangle&\lra&\cD\\
   \nu:\cB\langle X_1, X_2, \dots, X_N \rangle&\lra&\cB
   \end{eqnarray*}
    such that
  such that $\mu\circ\tau_{X_i}=\mu_i$ and  $\nu\circ\tau_{X_i}=\nu_{i}$ for all $i=1,\dots, N$ and the mixed moments of $\mu$ and $\nu$ are computed according to Definition \ref{def41}. The c-free additive convolution of  $\{(\mu_i, \nu_i)\}_{i=1}^N$ is the pair $(\mu_{c},\nu_{c})=\cfree_{i=1}^N (\mu_i, \nu_i)$, where
  \begin{eqnarray*}
  \nu_{c}&=&\nu\circ\tau_{X_1+X_2+\cdots+X_N}=\boxplus_{i=1}^n\nu_i\in\Sb\\
  \mu_{c}&=&\mu\circ\tau_{X_1+X_2+\cdots+X_N}:\cB\langle X_1+X_2+\dots+X_N\rangle\cong\bx\lra\cD
  \end{eqnarray*}

 \begin{defn}\label{def43}
  A pair $(\mu, \nu)\in\Sbd\times\Sb$ is said to be $\cfree$-infinite divisible if for any positive integer $N$ there exists $(\mu_N, \nu_N)\in\Sbd\times\Sb$ such that $(\mu, \nu)$ is the c-free additive  convolution of $N$ copies of $(\mu_N, \nu_N)$.
 \end{defn}

  \subsection{C-free cumulants and infinite divisibility}

 \begin{defn}\label{def43}
 The c-free cumulants of the pair $(\mu, \nu)\in\Sbd\times\Sb$ are the multilinear functions $\cka_{n, \mu,\nu}:\cB^n\lra\cD$ given by the recurrence:
  \begin{eqnarray*}
  \mu({\X}b_1{\X}b_2\cdots {\X}b_n)&=&
  \sum_{p=1}^n\sum_{\substack{1<j_1<\dots<j_{p}\\l_ p=n-1}}
  \mu({\X}b_1\cdots {\X}b_{l_1})\cdot\\
  &&\hspace{-2.6cm}\cka_{p, \mu,\nu}\left(b_{l_1+1}\nu({\X}b_{l_1+2}\cdots {\X}b_{l_2}),\dots, b_{l_{p-1}+1}\nu({\X}b_{l_{p-1}+2}\cdots {\X}b_{l_p}), b_n\right)\\
  \end{eqnarray*}
 \end{defn}

\begin{remark}
 As in the previous section, can reformulate Definition \ref{def43} in terms of free and c-free cumulants. More precisely, the pair  $(\mu, \nu)$ is  $\cfree$-infinitely divisible if for any positive integer $N$ there exists some $(\mu_N, \nu_N)\in\Sbd\times\Sb$ such that for all $m$ we have that $\displaystyle \kappa_{m,\nu}= N\kappa_{m, \nu_N}\ \text{and}\ \cka_{m, \mu,\nu}=N\cka_{m, \mu, \nu}$.
 
 \hfill $\square$
 %\begin{eqnarray*}
 %\kappa_{m,\nu}&=&N\kappa_{m, \nu_N}\\
 %\cka_{m, \mu,\nu}&=&N\cka_{m, \mu, \nu}
% \end{eqnarray*}
 \end{remark}

\noindent Define the map ${}^c\rho_{\mu,\nu}:\bx\lra\cD$ as the $\cB$-bimodule extension of 
\[ {}^c\rho_{\mu, \nu}({\X}b_1\cdots {\X}b_n)=\cka_{n, \mu,\nu}(b_1, \dots, b_n).\]

\begin{prop}\label{crhoprop}
 Suppose that $(\mu, \nu)\in\Sbd\times\Sb$ is $\cfree$-infinitely divisible. Then the restriction of ${}^c\rho_{\mu,\nu}$ to $\bx_0$ satisfies property (\ref{prop})(see Introduction).
\end{prop}

\begin{proof}

 Fix $N>1$ and suppose that $(\mu, \nu)$ is the c-free additive convolution of $n$ copies of $(\mu_N, \nu_N)$. As in the proof of Proposition \ref{rhoprop}, we will first show that
 \begin{equation}\label{eq43}
 \mu_N({\X}b_1\cdots {\X}b_n)=\frac{1}{N}\cka_{n,\mu,\nu}(b_1, \dots, b_n)+O(\frac{1}{N^2}).
 \end{equation}
 For $n=1$ the assertion is trivial. Suppose that (\ref{eq43}) holds true for $n<m$. Since the c-free cumulants of  are multilinear, for all $1=l_1<l_2<\dots<l_{p+1}<m$ and $Y_s=b_{l_s}\nu_N({\X}b_{l_s+1}\cdots {\X}b_{l_{s+1}-1})$, ($1\leq s \leq p$) we have that
  \[
    Y_s=\left \{
    \begin{array}{ll}
     b_{l_s}& \text{if}\ l_{s+1}=l_s+1\\
    O(\frac{1}{N})& \text{if}\ l_{s+1}=l_s+1
    \end{array}
    \right.
  \]
  hence
$\displaystyle
   \cka_{p,\mu_N, \nu_N}(Y_1, \dots, Y_p)
   =\frac{1}{N}\cka_{p, \mu, \nu}(Y_1, \dots, Y_p)
   =O(\frac{1}{N^2}),
   $\ 
    unless $l_{s+1}=l_s+1$ for all $s\in\{1, \dots, p\}$.

   Definition \ref{def43} gives
   \begin{eqnarray*}
    \mu_N({\X}b_1\cdots {\X}b_m)&=&\kappa_{m, \mu_N, \nu_N}(b_1, \dots, b_m)+\\ &&\hspace{.5cm}\sum_{s=1}^{m-1}\cka_{s, \mu_N, \nu_N}(b_1, \dots, b_s)\mu_N({\X}b_{s+1}\dots {\X}b_m)
    +O(\frac{1}{N^2})
    \end{eqnarray*}
    hence (\ref{eq43}) follows from the induction hypothesis. Therefore
    \begin{equation}\label{eq42}
    \lim_{N\lra\infty}N \mu_N({\X}b_1\cdots {\X}b_n)=\cka_{n, \mu, \nu}(b_1, \dots, b_n)
    \end{equation}

    Fix now a family $\{f_i(\X)\}_{i=1}^n$ in $\bx_0$. Then
    \begin{equation*}
    [{}^c\rho_{\mu,\nu}(f_j(\X)^\ast f_i(X))]_{i,j=1}^n
    =\lim_{N\lra\infty}[N\mu_N(f_j(\X)^\ast f_i(X))]_{i,j=1}^n=\geq 0.
    \end{equation*}
since each $\mu_N$ satisfies (\ref{prop}).
\end{proof}

\begin{lemma}\label{cfreeconstr}

  Let $\cB\subseteq\cD$ be unital inclusion of C$^\ast$-algebras, $\cK$ be a semi- inner-product $\cD$-bimodule and $\cH$ be a semi- inner-product $\cB$-bimodule. Consider
  \[
   \cE=\left(\cT(\cH)\otimes_\cB\cD\right)\oplus\left(1_\cB\otimes_\cB\Omega\cD\right)\oplus\left(\cT(\cH)\otimes_\cB\cK\right).
   \]

   Fix $\xi\in\cH$, $\eta\in\cK$ and $t\in\cL(\cH)$, $T\in\cL(\cK)$ selfadjoints. Define the maps $a_\xi$, $a^\ast_\xi$, $A_\eta$, $A^\ast_\eta$, $p(t)$, $P(T)$ from $\cL(\cE)$ given by:
   \begin{eqnarray*}
   &&\left\{
   \begin{array}{l}
   a_{\xi|1_\cB\otimes\Omega\cD}\equiv 0\\
   a_\xi(f_1\otimes\dots f_n)\otimes d=(\langle f_1, \xi\rangle f_2\otimes\dots\otimes f_n)\otimes d\\
   a_\xi(f_1\otimes\dots f_n)\otimes \eta=(\langle f_1, \xi\rangle f_2\otimes\dots\otimes f_n)\otimes \eta\\
  \end{array}
   \right.\\
      &&\left\{
   \begin{array}{l}
   a^\ast_{\xi|1_\cB\otimes\Omega\cD}\equiv 0\\
   a^\ast_\xi(f_1\otimes\dots f_n)\otimes d=(\langle f_1, \xi\rangle f_2\otimes\dots\otimes f_n)\otimes d\\
   a^\ast_\xi(f_1\otimes\dots f_n)\otimes \eta=(\langle f_1, \xi\rangle f_2\otimes\dots\otimes f_n)\otimes \eta\\
  \end{array}
   \right.\\
      &&\left\{
   \begin{array}{l}
   p(t)_{|1_\cB\otimes\Omega\cD}\equiv 0\\
   p(t)(f_1\otimes\dots f_n)\otimes d=(tf_1)\otimes f_2\otimes\dots\otimes f_n)\otimes d\\
   p(t)(f_1\otimes\dots f_n)\otimes \eta=(tf_1)\otimes f_2\otimes\dots\otimes f_n)\otimes \eta\\
  \end{array}
   \right.\\
    &&\left\{
   \begin{array}{l}
   A_\eta(1\otimes\zeta)=1\otimes\Omega\langle\zeta, \eta\rangle\\
   A_{\eta|\cE\ominus(1\otimes\cK)}\equiv 0\\
   \end{array}
   \right.
   \hspace{1.2cm}
   \left\{
   \begin{array}{l}
   A^\ast_\eta(1\otimes\Omega d)=1\otimes\eta d\\
   A^\ast_{\eta|\cE\ominus(1\otimes\Omega\cD)}\equiv 0\\
   \end{array}
   \right.\\
   &&\left\{
   \begin{array}{l}
   P(T)(1\otimes\zeta)=1\otimes T\zeta\\
   P(T)_{|\cE\ominus(1\otimes\cK)}\equiv 0\\
   \end{array}
   \right.\\
   \end{eqnarray*}
  Define also $I_1=\id_{(\cT(\cH)\otimes\cD)\oplus(\cT(\cH)\otimes\cK)}\oplus0\in\cL(\cE)$ and $I_2=\id_{1\otimes\Omega\cD}\oplus 0\in\cL(\cE)$.
 Consider $X=a_\xi+a_\xi^\ast+p(t)+\lambda_1I_1+A_\eta+A^\ast_\eta+\lambda_2I_2$ where $\lambda_1\in\cB$, $\lambda_2\in\cD$ selfadjoint.
 Then the free and c-free cumulants of $X$ with respect to $(\theta, \varphi)=\left(\langle \cdot1\otimes\Omega, 1\otimes\Omega\rangle,
 \langle \cdot 1\otimes 1, 1\otimes 1\rangle\right)$ are given by the following relations:
 \begin{eqnarray*}
 &&\left\{
 \begin{array}{l}
 \kappa_{1, X}(b_1)=\lambda_1b_1\\
 \kappa_{n,X}(b_1, \dots, b_n)=
 \langle a_\xi b_1 p(t)b_2\cdots p(t)b_{n-1}a^\ast_\xi b_1(1\otimes 1),
 1\otimes 1\rangle
 \end{array}
 \right.\\
 &&\left\{
 \begin{array}{l}
 {}^c\kappa_{1, X}(b_1)=\lambda_2b_1\\
 {}^c\kappa_{n,X}(b_1, \dots, b_n)=
 \langle A_\eta b_1 P(T)b_2\cdots P(T)b_{n-1}A^\ast_\eta b_1(1\otimes \Omega),
 1\otimes \Omega\rangle
 \end{array}
 \right.
 \end{eqnarray*}
\end{lemma}

\begin{proof}
 The results are trivial for $\kappa_{1, X}$ and ${}^c\kappa_{1, X}$. For $\kappa_{n, X}$ note that $a_\xi, a^\ast_\xi, p(t)$ map $\cT(\cH)\otimes 1$ in $\cT(\cH)\otimes 1$ and, since $1\otimes1\in\cT(\cH)\otimes 1$, the result reduces to Lemma \ref{}.

 To prove the formula for ${}^c\kappa_{n, X}$, let us first note $V_0=a_\xi+a^\ast_\xi+p(t)+\lambda_1I_1$, $V_1=A^\ast_\eta$, $V_2=A_\eta$, $V_3=P(T)$, $V_4=\lambda_2I_2$ and
 \[
 J(n)=\{\vu=(u_n, \dots, u_1): 0\leq u_k\leq 4\}.
 \]
 Finally, for $b_1, \dots, b_n\in\cB$, denote $\epsilon(\vu)=\theta(V_{u_n}b_n\cdots V_{u_1}b_1).$

 Note that  $V_{k|\cT(\cH)\otimes 1}\equiv0$ ($1\leq k\leq 4$) and $V_{0|1\otimes\Omega\cD}\equiv0$. The latest implies that $\epsilon(\vu)=0$ unless $V_{u_1}\in\{A_\eta^\ast, \lambda_2I_2\}$, hence
 \begin{eqnarray*}
 \theta(Xb_n\cdots Xb_1)&=&\sum_{\vu\in J(n)}\epsilon(\vu)=
 \sum_{\substack{\vu\in J(n)\\ u_1=1}}\epsilon(\vu)+
 \sum_{\substack{\vu\in J(n)\\ u_1=4}}\epsilon(\vu)
 \\
 &=&
 \sum_{\substack{\vu\in J(n)\\ u_1=4}}\epsilon(\vu) +\sum_{\vu\in J(n-1)}\theta(V_{u_{n-1}}b_n\cdots V_{u_1}b_2)\lambda_2
 \end{eqnarray*}

 Suppose that $V_{u_1}=A^\ast_\eta$; then $V_{u_1}b_1(1\otimes\Omega)=1\otimes\eta b_1$, hence $\epsilon(\vu)$ cancels unless $V_{u_2}\in\{ V_0, A_\eta, P(T)\}$. Let $p=\text{min}\{s: s>1, u_s\neq0\}$. It follows that $V_{u_p}\in\{A_\eta, P(T)\}$. Since the restrictions of $A_\eta, P(T)$ to $\cT(\cH)\otimes\cK$ are $0$, we have
 \begin{eqnarray*}
 V_{u_p}b_p\cdots V_{u_1}b_1(1\otimes\Omega)&=&V_{u_p}b_p\left(V_0b_{p-1}\cdots V_0b_2\right)V_{u_1}b_1(1\otimes\Omega)\\
 &=&
 V_{u_p}b_p\varphi(Xb_{p-1}\cdots Xb_2)V_{u_1}b_1(1\otimes\Omega).
 \end{eqnarray*}

If $V_{u_p}=P(T)$ and $s=\text{min}\{q: p<q\leq n, u_q\neq 0\}$, from a similar argument as above we have that $u_{p+1}=\dots=u_{q-1}=0$ and
\begin{eqnarray*}
V_{u_q}b_q\cdots V_{u_1}b_1(1\otimes\Omega)&=& \\
&&\hspace{-2cm}P(T)b_q\varphi(Xb_{q-1}\cdots Xb_{p+1})P(T)b_p\varphi(Xb_{p-1}\cdots Xb_2)A^\ast_\eta b_1(1\otimes\Omega)
\end{eqnarray*}

 Note also that, for all $b_1, \dots, b_m\in \cB$, one has
 \begin{eqnarray*}
 Tb_m\cdots  Tb_2A^\ast_\eta b_1(1\otimes\Omega)&\in &1\otimes \cK\\
 A_\eta b_m T b_{m-1}\cdots T b_{2} A^\ast_\eta b_1(1\otimes\Omega)&\in &1\otimes \Omega\cD,
 \end{eqnarray*}
 hence $\epsilon(\vu)$ cancels unless there exists some $j>1$ such that $V_{u_j}=A_\eta$.

 Using the results above, we have that
 \begin{eqnarray*}
 \sum_{\substack{\vu\in J(n)\\u_1=4}}\epsilon(\vu)
 &=&
 \sum_{s=1}^n
 \sum_{\substack{p_1\cdots<p_s\\ 1<p_1, p_s\leq n}}
 \sum_{\substack{u\in J(n)\\ u_1=4, u_{p_s}=2\\ u_{p_l}=1, l\neq s}}
 \epsilon(\vu)\\
 &&
 \hspace{-1cm}=
 \sum_{s=1}^n
 \sum_{\substack{p_1\cdots<p_s\\ 1<p_1, p_s\leq n}}
 \theta(Xb_n\cdots Xb_{p_{s}+1})\cdot \theta(A_\eta b_{p_s}\varphi(Xb_{p_s-1}\cdots Xb_{s-1}+1)\cdot\\
 &&\hspace{-1cm}P(T)b_{p_{s-1}}\varphi(Xb_{p_{s-1}-1}\cdots X b_{p_{s-2}+1})\cdots
 P(T)b_{p_1}\varphi(Xb_{p_1}-1\cdots Xb_2)A^\ast_\eta b_1)
 \end{eqnarray*}

  Comparing these relations with the definition of the c-free cumulants, we have q.e.d..
\end{proof}

\begin{thm}\label{thm45}
 The pair $(\mu, \nu)\in\Sbd\times\Sb$ is $\cfree$-infinitely divisible if and only if $\nu$ is $\boxplus$-infinitely divisible and the restriction of
  ${}^c\rho_{\mu, \nu}$ to $\bx_0\rangle$ satisfies property (\ref{prop}).
\end{thm}
\begin{proof}

 Suppose that $\nu$ is infinitely divisible (hence, from Theorem \ref{rhothm}, the restriction of $\rho_\nu$ to $\bx_0$ is positive) and that the restriction of ${}^c\rho_{\mu, \nu}$ to $\bx_0$ satisfies (\ref{prop}).

 Let $\cH$ be the left $\cB$-bimodule $\bx_0$ with the pairing
 \[\langle f(\X), g(\X)\rangle_\cH=\rho_\nu(g(\X)^\ast f(\X))\]
 and $\cK$ be the left $\cD$-module $\bx_0\otimes_\cB\cD$ with the pairing
 \[
 \langle f(\X)\otimes d_1, g(\X)\otimes d_2\rangle_\cD=d_2^\ast\cdot {}^c\rho_{\mu, \nu}(g^\ast(\X)f(\X))d_1.
 \]
  Since the maps $\rho_{\nu}$ and ${}^c\rho_{\mu, \nu}$ staify (\ref{prop}), we have that $\cH$ and $\cK$ are semi-inner-product $\cB$-, respectively $\cD$-modules.

  Define
  \begin{eqnarray*}
  t\in\cL(\cH), & \ & t(f(\X))=\X f({\X})\hspace {3cm} \\
  T\in\cL(\cK), & \ &  T(f({\X})\otimes d)={\X}f({\X})\otimes d \hspace{3cm}
 \end{eqnarray*}
 and note that $t$ and $T$ are selfadjoint maps. Also let $\lambda_1=\rho_\nu({\X})\in\cB$ and\\ $\lambda_2={}^c\rho_{\mu. \nu}({\X})\in\cD$.

 As in Lemma \ref{cfreeconstr}, consider
   \[
   \cE=\left(\cT(\cH)\otimes_\cB\cD\right)\oplus\left(1_\cB\otimes_\cB\Omega\cD\right)\oplus\left(\cT(\cH)\otimes_\cB\cK\right)
   \]
   and $V\in\cL(\cE)$ given by
   \[V=a_{\X}+a_{\X}^\ast+p(t)+\lambda_1I_1+A_{{\X}\otimes 1}+A^\ast_{{\X}\otimes1}+\lambda_2I_2. \]

   Denoting $\varphi(\cdot)=\langle \cdot (1\otimes1), 1\otimes 1\rangle$, $\theta(\cdot)=\langle \cdot (1\otimes\Omega), 1\otimes\Omega\rangle$ and applying Lemma \ref{cfreeconstr} we have that the free and c-free cumulants of $V$ with respect to $(\theta, \varphi)$ are given by:
  \begin{eqnarray*}
  \kappa_{1, {\X}}(b_1)&=&\lambda_1b_1=\rho_{\nu}({\X})\\
  &=&
  \kappa_{1,\nu}(b_1)\\
  \kappa_{n, {\X}}(b_1, \dots, b_n)
  &=&
  \langle a_{\X} b_1 p(T)b_2\cdots p(t)b_{n-1}a^\ast_{\X}b_n (1\otimes 1), 1\otimes 1\rangle\\
  &=&
  \langle a_{\X} (b_1{\X}b_2\cdots {\X}b_{n-1} {\X}b_n \otimes 1), 1\otimes 1\rangle\\
  &=&
  \langle\rho_\nu({\X}b_1\cdots {\X}b_n)\otimes1, 1\otimes1\rangle\\
  &=&
  \kappa_{n, \nu}(b_1, \dots, b_n),
  \end{eqnarray*}
 respectively  by
 \begin{eqnarray*}
 \cka_{1, V}(b_1)
 &=&
 \lambda_2b_1={}^c\rho_{\mu, \nu}({\X}b_1)\\
 &=&\cka_{1, \mu,\nu}(b_1)\\
 \cka_{n, V}(b_1, \dots, b_n)
 &=&
 \langle
 A_{{\X}\otimes1}b_1P(T)b_2\cdots P(T)b_{n-1}A^\ast_{{\X}\otimes1}b_1(1\otimes\Omega),
 1\otimes\Omega
 \rangle\\
 &=&
 \langle
 A_{{\X}\otimes1}
 (b_1{\X}b_2\cdots {\X}b_n\otimes1),
 1\otimes\Omega
 \rangle\\
 &=&
 \langle
 {}^c\rho_{\mu,\nu}({\X}b_1\cdots {\X}b_n)\otimes\Omega,
 1\otimes\Omega
 \rangle\\
 &=&
 {}^c\rho_{\mu,\nu}({\X}b_1\cdots {\X}b_n)\\
 &=&
 \cka_{n, \mu,\nu}(b_1, \dots, b_n).
 \end{eqnarray*}

 Fix $N>0$ and define
 \[
V_N=\frac{1}{\sqrt{N}}a_{\X}+\frac{1}{\sqrt{N}}a^\ast_{\X}+p(t)+\frac{1}{N}\lambda_1I_1+\frac{1}{\sqrt{N}}A_{{\X}\otimes1}+
\frac{1}{\sqrt{N}}A^\ast_{{\X}\otimes1}
+P(T)+\frac{1}{N}\lambda_2I_2.
\]
Using again Lemma \ref{cfreeconstr}, similar computations as above give ($n\geq1$):
\begin{eqnarray*}
\kappa_{n, V_N}
&=&
\frac{1}{N}\kappa_{n, V}\\
\cka_{n, V_N}
&=&
\frac{1}{N}
\cka_{n,V}
\end{eqnarray*}
hence q.e.d..

The converse implication is Proposition \ref{crhoprop}.
\end{proof}

\section{The non-commutative $R-$ and ${}^cR$-transforms and non-commutative free Levy-Hincin formula}

 \subsection{The $R$- and ${}^cR-$transforms and free infinite divisibility: scalar case}
 \begin{defn}
 Let $(\mu, \nu)$ be a pair of compactly supported measures on $\mathbb{R}$. If $M_\mu(z)$, $M_\nu(z)$ are the moment-generating series for $\mu$, respectively $\nu$, that is
 \begin{eqnarray*}
 M_\mu(z)
 &=&
 \sum_{n=0}^\infty \int_{\mathbb{R}}t^{n}d\mu(t)z^{n}\\
 M_\nu(z)
 &=&
 \sum_{n=0}^\infty \int_{\mathbb{R}}t^{n}d\nu(t)z^{n}
 \end{eqnarray*}
 then the $R$-transform of $\nu$, respectively the ${}^cR$-transform of $(\mu,\nu)$  are the  analytic functions $R_\nu$, ${}^cR_{\mu,\nu}$ given by
 \begin{eqnarray}
 M_\nu(z)-1
 &=&
 R_\nu\left(z M_\nu(z)\right)\label{Rtr}\\
 (M_\mu(z)-1)M_\nu(z)
 &=&
 M_\mu(z)\cdot{}^cR_{\mu,\nu}\left(zM_\nu\right(z))\label{cRtr}.
 \end{eqnarray}

\end{defn}

 If $\textit{A}$ is a $\ast$-algebra and $\varphi:\textit{A}\lra\mathbb{C}$ is a positive linear functional, then a selfadjoint element $a\in\textit{A}$ determines a compactly supported real measure $\mu_a$ via
 \[\int_{\mathbb{R}}t^nd\mu_a(t)=\varphi(a^n).\]
 For convenience, we will also use the notation $R_a$ for $R_{\mu_a}$, respectively the notation ${}^cR_a$ for ${}^cR_{\mu_a, \nu_a}$ for the case that the the $\ast$-algebra $\textit{A}$ endowed with two positive linear functionals.

  The key property of the $R$- and ${}^cR$-transforms is the linearization of the free convolution (\cite{vdn}, \cite{speichernica}), respectively c-free additive convolution: if $a$ and $b$ are free, respectively c-free, selfadjoint elements from $\textit{A}$, then
  \begin{eqnarray}
  R_{a+b}(z)
  &=&R_{a}(z)+R_b(z)\label{Rtransf}\\
  {}^cR_{a+b}
  &=&{}^cR_{a}(z)+{}^cR_b(z)\label{cRtransf}
  \end{eqnarray}

   The goal of this material is to give a non-commutative analogue for the following theorem (\cite{speichernica}, Theorem 13.16 and  \cite{Krystek}, Theorem 5):
   \begin{thm}
    Let $\nu$ be a compactly supported supported probability measure on $\mathbb{R}$ and let $R_\nu(z)=\sum_{n=1}^\infty \kappa_{n}z^n$ be the Taylor expansion of its $R$-transform. Then the following statements are equivalent:
    \begin{enumerate}
    \item[(1)] $\nu$ is infinitely divisible.
    \item[(2)] The sequence $\{\kappa_n\}_n\geq2$ is positive definite, i.e. there exists some real measure $\sigma$ such that
    \[\kappa_n=\int_\mathbb{R}t^{n-2}d\sigma\]
    \item[(3)] The $R$-transform of $\nu$ is of the form
    \[\frac{1}{z}R_\nu(z)=\kappa_1+\int_{\mathbb{R}}\frac{z}{1-tz}d\rho(t),\]
    for some finite measure $\rho$ on $\mathbb{R}$ with compact support.
    \end{enumerate}
    Moreover, if $(\mu, \nu)$ is a pair of compactly supported measures on $\mathbb{R}$, then $(\mu, \nu)$ is c-free infinitely divisible if and only if $\nu$ is infinitely divisible and ${}^cR_{\mu,\nu}$ satisfies the condition (3) from above.
   \end{thm}

  \subsection{Non-Commutative functions}\label{ncfcn}
%subsection rewritten by VV
For stating our main result, we introduce the language of noncommutative \cite{ncfound} or fully matricial \cite{voi09}
functions, see also the pioneering work \cite{T72b,T73}.
For a vector space $\cV$ over $\mathbb{C}$, we let
$\cV^{n\times m} = \cV \otimes_{\mathbb{C}} M_{n\times m}(\mathbb{C})$ denote $n \times m$ matrices over $\cV$ (in literature - for example in \cite{ncfound} - on $\mat{\cV}{n}$ is used the algebra structure induced by the tensor algebra $\ten{\cV}$ over $\cV$; in order to avoid confusion, when $\cV$ is an algebra, we will use the notation $M_n(\cV)$ for the \emph{algebra} $M_n(\cV)$ of $n\times n$ matrices over $\cV$).
We define the {\em noncommutative space} over $\cV$ by
$\ncspace{\cV} = \displaystyle\coprod_{n=1}^\infty \mat{\cV}{n}$.
We call $\Omega \subseteq \ncspace{\cV}$ a {\em noncommutative set} if it is closed under direct sums.
Explicitly, denoting $\Omega_n = \Omega \cap \mat{\cV}{n}$,
we have
\[
a \oplus b =\begin{bmatrix}
a & 0\\
0 & b\end{bmatrix}\in \Omega_{n+m}
\]
for all $a \in \Omega_n$, $b \in \Omega_m$.
Notice that matrices over $\mathbb{C}$ act from the right and from the
left on matrices over $\cV$  by the standard rules of
matrix multiplication.

A noncommutative set $\Omega\subseteq\ncspace{\cV}$ is called {\em upper admissible} if for all
$a \in \Omega_n$, $b \in \Omega_m$ and all $c \in \rmat{\cV}{n}{m}$,
there exists $\lambda \in \mathbb{C}$, $\lambda \neq 0$, such that
\[
\begin{bmatrix} a & \lambda c \\ 0 & b \end{bmatrix} \in \Omega_{n+m}.
\]
This notion is crucial since it is used to define the (right) noncommutative difference-differential operators
by evaluating a noncommutative function on block upper triangular matrices.
We will encounter only the following upper admissible noncommutative sets:
\begin{enumerate}
\item
The set $\Nilp(\cV) = \coprod_{n=1}^\infty \Nilp(\cV;n)$ of nilpotent matrices over $\cV$.
Here the set $\Nilp(\cV;n)$ of nilpotent $n \times n$ matrices over $\cV$ consists of
all $a \in \mat{\cV}{n}$ such that $a^r = 0$ for some $r$, where we view $a$ as a matrix over
the tensor algebra $\ten{\cV}$ of $\cV$ over $\mathbb{C}$.
This is equivalent to $t a t^{-1}$ being strictly upper triangular for some $t \in M_n(\mathbb{C})$
(the equivalence follows from Engel's Theorem --- notice that we can restrict ourselves to
the finite dimensional subspace of $\cV$ spanned by the elements of $a$).
\item
A noncommutative ball $\ball(\cA,\rho) = \left\{a \in \ncspace{\cA} \colon \|a\| < \rho\right\}$ centered in zero and
of radius $\rho>0$
over a C$^\ast$-algebra $\cA$ ($\cA$ could have been replaced by any operator space with the corresponding
operator space norm).
\item
The upper and lower fully matricial half-planes $\mathbb{H}^{+}(\ncspace{A})$ and $\mathbb{H}^{-}(\ncspace{A})$ over a C$^\ast$-algebra $A$, where if $\mathcal{C}$ is a C$^\ast$-algebra, then
\begin{eqnarray*}
\mathbb{H}^{+}(\mathcal{C})&=&\{a\in\mathcal{C}, \Im{a}=\frac{(a-a^\ast)}{2}>0 \}\\
\mathbb{H}^{-}(\mathcal{C})&=&\{a\in\mathcal{C}, \Im{a}=\frac{(a-a^\ast)}{2}<0\}
\end{eqnarray*}
and $\mathbb{H}^{\pm}(\ncspace{A})=\coprod_{n=1}^\infty\mathbb{H}^{\pm}(M_n(A))$.
\end{enumerate}

Let $\cV$ and $\cW$ be vector spaces over $\mathbb{C}$, and let
$\Omega\subseteq\ncspace{\cV}$ be a noncommutative set. A
mapping $f \colon \Omega \to \ncspace{\cW}$ with
$f(\Omega_n) \subseteq \mat{\cW}{n}$ is called a
\emph{noncommutative function} if $f$ satisfies the following two
conditions:
\begin{itemize}
\item $f$ \emph{respects direct sums}:
%\begin{equation} \label{eq:dirsums}
$f(a \oplus b) = f(a) \oplus f(b)$
%\end{equation}
for all $a,b \in\Omega$.
\item $f$ \emph{respects similarities}:
if $a \in \Omega_n$ and $s \in \mat{\mathbb{C}}{n}$ is invertible with
$sas^{-1} \in \Omega_n$, then
%\begin{equation} \label{eq:simsim}
$f(sas^{-1}) = s f(a) s^{-1}$.
%\end{equation}
\end{itemize}
We will denote $f_n = f|_{\Omega_n} \colon \Omega_n \to \mat{\cW}{n}$.
While we will not need this fact, it is important to notice that
the two conditions in the definition of a noncommutative function can be
actually replaced by a single one:
a mapping $f \colon \Omega \to \ncspace{\cW}$ with
$f(\Omega_n) \subseteq \mat{\cW}{n}$ respects direct sums
and similarities if and only if it \emph{respects intertwinings}:
for any $a\in\Omega_n$, $b\in\Omega_m$, and $t\in
\rmat{\mathbb{C}}{n}{m}$ such that $at=tb$, one has $f(a)t=tf(b)$.
This last condition goes back to \cite{T72b,T73}.

Let $\alpha \colon \cV^k \longrightarrow \cW$  be a multilinear mapping;
we set
\[
\talpha = \talpha_n \overset{\text{def}}{=} \alpha \otimes \id_n
\colon \mat{\left(\cV^{\otimes k}\right)}{n} \longrightarrow \mat{\cW}{n}
\]
for each $n \in \mathbb{N}$.
It is clear that the mapping $a \mapsto \talpha(a^k)$,
where we view $a$ as a matrix over $\ten{\cV}$,
respects direct sums and similarities,
so that it defines a noncommutative function from $\ncspace{\cV}$ to $\ncspace{\cW}$.
Therefore for any linear mapping $\phi \colon \ten{\cV} \longrightarrow \cW$, we obtain a noncommutative
function defined by
\begin{equation} \label{eq:ncfunfromforms}
f(a) = \tphi \left((\mathds{1} - a)^{-1}\right)
         = \sum_{k=0}^\infty \tphi(a^k)
\end{equation}
(where the notation $\mathds{1}$ should be understood as $\id_n$ componentwise, i.e. $f(a)=\tphi_n \left((\id_n - a)^{-1}\right)$ for $a \in \mat{\cV}{n}$),
except that we have to make sense of the infinite sum on the right-hand side.
\begin{enumerate}
\item
If $a \in \Nilp(\cV)$ then the sum is finite, so that $f$ is always a noncommutative function
on $\Nilp(\cV)$.
\item
If $\phi \colon \ten{\cA} \longrightarrow \cC$, where $\cA$ and $\cC$ are $C^*$ algebras, and we have
an exponential growth estimate: $\|\phi|_{\cA^{\otimes k} \to \cC}\|_{\text{cb}} \leq \alpha \beta^k$
(where $\cA^{\otimes k}$ is considered with the Haagerup tensor norm (see \cite{paulsen}, Chapter 17) and $\|\cdot\|_{\text{cb}}$ denotes the completely bounded norm),
then the series defining $f$ converges absolutely and uniformly on any noncommutative ball
$\ball(\cA,r)$ over $\cA$
of radius $r < 1/\beta$,
so that $f$ is a noncommutative function on the noncommutative ball $\ball(\cA,1/\beta)$.
\end{enumerate}

There is --- in a sense --- a converse to this construction that we briefly describe,
though we will make no real use of it here. A noncommutative function $f$
admits a series expansion
\begin{equation} \label{eq:tt}
f(a) = \sum_{k=0}^\infty \widetilde{\Delta^k_R f(\underset{k+1}{\underbrace{0,\ldots,0}})}(a^{\otimes k}).
\end{equation}
More precisely:
\begin{enumerate}
\item
If $f$ is a noncommutative function on $\Nilp(\cV)$, then the sum is finite and the equality holds
everywhere.
\item
If $f$ is a noncommutative function on a noncommutative ball $\ball(\cA,\rho)$ over a $C^*$ algebra
$\cA$ with values in a noncommutative space $\ncspace{\cC}$ over a $C^*$ algebra $\cC$,
which is bounded on noncommutative balls of radius less then $\rho$, then the series
converges to $f$, uniformly on every noncommutative ball of radius less than $\rho$.
(The convergence still holds if $f$ is only assumed to be locally bounded in every matrix dimension
separately, but then it is no longer uniform across matrix dimensions.)
\end{enumerate}
The multilinear forms
$\Delta^k_R f(\underset{k+1}{\underbrace{0,\ldots,0}}) \colon \cV^k \longrightarrow \cW$ are the values
at $(0,\ldots,0)$ of the $k$th order noncommutative difference-differential operators applied to $f$.
They are uniquely determined, and can be calculated directly by evaluating $f$ on upper triangular matrices:
\begin{multline*}
f\left(
\left[
 \begin{array}{ccccc}
 0&a_1&0&\dots&0\\
 0&0&a_2&\dots&0\\
 \dots&\dots&\dots&\dots&\dots\\
 0&0&0&\dots&a_k\\
 0&0&0&\dots&0
\end{array}
\right]
%\begin{bmatrix}
%0 & a_1 & 0 & \cdots & 0\\
%0   & 0 & \ddots & \ddots & \vdots\\
%\vdots & \ddots & \ddots & \ddots & 0 \\
%\vdots &  & \ddots & 0 & a_k\\
% 0 & \cdots & \cdots & 0 & 0
%\end{bmatrix}
\right)\\
=\begin{bmatrix} f(0) & \Delta_Rf(0,0)(a_1) & \cdots &
\cdots &
\Delta_R^k f(0,\ldots,0)(a_0,\ldots,a_k)\\
0   & f(0) & \ddots & & \Delta_R^{k-1} f(0,\ldots,0)(a_2,\ldots,a_k)\\
\vdots & \ddots & \ddots & \ddots & \vdots \\
\vdots &  & \ddots &  f(0) & \Delta_Rf(0,0)(a_k)\\
 0 & \cdots & \cdots &  0 & f(0)
\end{bmatrix}.
\end{multline*}

\subsection{The non-commutative $B$-, $R$- and ${}^cR$-transforms}${}$\\

In particular, we can apply the construction \eqref{eq:ncfunfromforms}
to $\mu\in\Sbd$, and define
\begin{equation*}
 M_\mu(b)=\tmu\left((\mathds{1}-\X b)^{-1}\right)=\sum_{k=0}^\infty \tmu\left((\X b)^{n}\right),
\end{equation*}
where the notation $\mathds{1}$ is the one from the previous section and, as before, we identify the monomials $(\X b)^n$ from $\bx_0$ to $b^{\otimes n}$ from the tensor algebra $\ten{\cB}$.
$M_\mu$ is a noncommutative function with values in $\ncspace{\cD}$.
It is always defined on $\Nilp(\cB)$.
If $\cD$ is a C$^\ast$-algebra and $\mu\in\Sbd^0$,
that is there is some $C>0$ such that $\|\mu_k\|_{\text{cb}} \leq C^{k+1}$,
where $\mu_k$ is the restriction of $\mu$ to the subspace of $\bx$ spanned by words with exactly $k$
occurrences of the symbol $X$ (that can be identified with $\cB^{\otimes k}$),
then $M_\mu$ is also defined on a noncommutative ball in $\cB$ with sufficiently small radius. Remark that \eqref{eq:ncfunfromforms} also implies  $\Delta_R^k M_\mu(0,\ldots,0) = \mu_k$.

Consider now $\nu\in\Sb$ and $\mu\in\Sbd$. We define the linear maps $\rho_\nu$, ${}^c\rho_{\mu,\nu}$ as in Section \ref{section3}, respectively Section \ref{cfreeinf}, that is the $\cB$-bimodule extensions of
\[
\rho_\nu(\X b_1\cdots \X b_p)=\kappa_{\nu,p}(b_1, \dots, b_p),\ \hspace{.7cm}
{}^c\rho_{\mu,\nu}(\X b_1\cdots \X b_p)={}^c\kappa_{\mu,\nu,p}(b_1, \dots, b_p)
\]
and the map $\beta_\mu$ as the $\cB$-bimodule extension of
\[
\beta_\mu(\X b_1\cdots \X b_p)= B_{\mu,p}(b_1,\cdot, b_p).
\]
Using again \eqref{eq:ncfunfromforms}, we define the noncommutative functions $R_\nu$ (with values in $\ncspace{\cB}$)
and  $B_\nu$, ${}^cR_{\mu,\nu}$ (with values in $\ncspace{\cD}$) via the equations
\begin{eqnarray*}
R_\nu(b)&=&\widetilde{\rho_\nu}\left((\mathds{1}-\X b)^{-1}\right)-\mathds{1}\\
{}^cR_{\mu,\nu}(b)&=&\widetilde{{}^c\rho_{\mu,\nu}}\left((\mathds{1}-\X b)^{-1}\right)-\mathds{1}\\
B_\mu(b)&=&\widetilde{\beta_\mu}\left((\mathds{1}-\X b)^{-1}\right)-\mathds{1}.
\end{eqnarray*}
Remark again that \eqref{eq:ncfunfromforms} implies $\Delta^p_R\cR_\nu(0, \dots, 0)=\kappa_{
 \nu,p}$, $\Delta^p_R{}^c\cR_{\mu,\nu}(0,\dots, 0)=\cka_{\mu,\nu,p}$ and $\Delta^p_R B_\mu(0,\dots,0)=B_{\mu,p}$.

 \begin{prop}\label{tensorprop}
 For all $(\mu,\nu)\in\Sbd\times\Sb$ we have that
 \begin{eqnarray*}
 \id_{n}\otimes \rho_\nu&=&\rho_{\nu^{(n)}}\\
 \id_{n}\otimes {}^c\rho_{\mu,\nu}&=&\rho_{\mu^{(n)},\nu^{(n)}}\\
 \id_{n}\otimes \beta_{\mu}&=&\beta_{\mu^{(n)}}.
 \end{eqnarray*}
 \end{prop}
\begin{proof}
 For the first property, it suffices to show that for all positive integers  $m$, all $1\leq i_k, j_k\leq n$  ($k=1,\dots, n$) and all $b_1,\dots, b_m\in\cB$ we have that
\[\id_n\otimes\rho_\nu(\X B_1\cdots \X B_m)=\rho_{\nu^{(n)}}(\X B_1\cdots \X B_m),\]
 where $B_k=e_{i_k, j_k}\otimes b_k$, for $e_{i_k, j_k}$ the complex $n\times n$ matrix with the $(i_k, j_k)$ entry 1 and all others 0.

For $m=1$ the property is trivial. Suppose now that the assertion holds true for all $m<N$. The definition of free cumulants gives
\begin{eqnarray*}
\rho_{\nu^{(n)}}(\X B_1\cdots \X B_N)
&=&
\nu^{(n)}(\X B_1\cdots \X B_N)-\\
&&\hspace{-3.5cm}
\sum_{p=1}^{N-1}\sum_{{s_1<\cdots<s_p\atop s_1=1, s_p\leq N}}\rho_{\nu^{(n)}}
\left(
\X B_1 \nu^{(n)}(\X B_2\cdots \X B_{s_2-1})\cdots \X B_{s_p}\nu^{(n)}(\X B_{s_p+1}\cdots \X B_N)
\right).
\end{eqnarray*}
 From the induction hypothesis, the right-hand side cancels unless $i_k=J_{k-1}$ for all $1<k\leq N$; in this case, the equation above becomes
\begin{eqnarray*}
\rho_{\nu^{(n)}}(\X B_1\cdots \X B_N)
&=&
e_{i_1,j_N}\otimes\nu(\X b_1\cdots \X b_N)-\\
&&\hspace{-3.5cm}
\sum_{p=1}^{N-1}\sum_{{s_1<\cdots<s_p\atop s_1=1, s_p\leq N}}\rho_{\nu^{(n)}}
\left(
e_{i_1,j_N}\otimes[
\X b_1 \nu(\X b_2\cdots \X b_{s_2-1})\cdots \X b_{s_p}\nu(\X b_{s_p+1}\cdots \X b_N)]
\right)\\
&=&
e_{i_1,j_N}\otimes\rho_\nu(\X b_1\cdots \X b_N)\\
&=&
\id_n\otimes\rho_\nu(\X B_1\cdots \X B_N)
\end{eqnarray*}
hence the conclusion.
 %%%%%%%%%%%%%%%%%%%%%%%%%%%%

 %%%%%%%%%%%%%%%%%%%%%%%%%%%%%%%%%%%%%%%%%%%%%%%%%%%%%%%%

\noindent  The argument for $^c\rho_{\mu,\nu}$ and $\beta_{\mu}$ is analogous.
\end{proof}

Proposition \ref{tensorprop} implies that
\[
R_\nu^{(n)}(b)
      =
      \sum_{n=1}^\infty\kappa_{n, \nu^{(n)}}(b,\dots, b)\]
and the analogous relations for the components of $\CR_{\mu,\nu}$ and $B_\mu$. From the moment-cumulant recursions, we have then the following
\begin{cor}
The non-commutative functions $R_\nu$,  $B_\nu$ and ${}^cR_{\mu,\nu}$ satisfy the equations:
  \begin{eqnarray}
    M_\mu(b)-\mathds{1}
 &=&B_\mu(b)\cdot M_\mu(b)\label{eq400}\\
 M_\nu(b)-\mathds{1}
 &=&
 R_\nu\left(b M_\nu(b)\right)\label{eq41}\\
 (M_\mu(b)-\mathds{1})\cdot M_\nu(b)
 &=&
 M_\mu(b)\cdot{}^cR_{\mu,\nu}\left(b M_\nu(b)\right)\label{eq42}.
 \end{eqnarray}
\end{cor}

\begin{remark}\label{remark55}
 If $(\mu,\nu)\in\Sbd^0\times\Sb^0$ then $R_\nu$, $\CR_{\mu,\nu}$ and $B_\mu$ are  well-defined on a non-commutative ball (as described in Section \ref{ncfcn}) centered in $0$.

\end{remark}

\begin{proof}
  For $B_\mu$, the assertion is trivial since $B_\mu(b)$ is well-defined if $M_\mu$ is invertible. For $R_\nu$ and  $\CR_{\mu,\nu}$ we will use combinatorial techniques, namely the Moebius inversion formula for the partially ordered set of non-crossing partitions. First, we need the following general property (the Moebius inversion formula).

  \textbf{Proposition:} \textit{Let $P$ be a finite partially ordered set and $K$ a complex vector space. Then there exist a map}
  $\text{moeb}:P\times P\lra\mathbb{R}$
  \textit{such that if the maps $f,g:P\lra K$ have the property}
  \[\displaystyle f(\pi)=\sum_{\substack{\sigma\in P\\ \sigma\leq \pi}}g(\sigma), \pi\in P\]
  \textit{then} \[\displaystyle g(\pi)= \sum_{\substack{\sigma\in P\\ \sigma\leq \pi}}f(\sigma)\cdot\text{moeb}(\sigma, \pi).\]

  We will apply the above result on the partially ordered set $\mathcal{NC}$ of non-crossing partitions. By a partition on the ordered set $\langle n\rangle=\{1,2,\dots, n\}$  we will understand a collection of mutually disjoint subsets of $\langle n \rangle$, $\gamma=(B_1,\dots, B_q)$ called {\em blocks} whose union is the entire set $\langle n \rangle$. A crossing is a sequence $i < j < k < l$
from $\langle n \rangle$ with the property that there exist two different blocks $B_r$ and $B_s$ such that
$i, k\in B_r$ and $j, l\in B_s$. A partition that has no crossings will be called non-crossing.
The set of all non-crossing partitions on $\langle n \rangle$ will be denoted by NC(n). For $\gamma\in NC(n)$ a block $B = (i_1, . . . , i_k)$ of $\gamma$ will be called {\em interior} if there exists
another block $D\in\gamma$ and $i, j\in D$ such that $i < i_1, i_2,\dots, i_k < j$. A block will be
called {\em exterior} if is not interior. Each $NC(n)$ has a lattice structure with respect to block refinement with biggest element $1_m$- the partition with a single block; the corespondent Moebius function satisfies $|\text{moeb}(\pi, 1_n)|\leq 4^n$ (see \cite{speichernica}, Lecture 13). Define $\mathcal{NC}=\coprod_{n=1}^\infty NC(n)$.

 Finally, let $(\mu,\nu)\in\Sbd\times\Sb$, $b\in\cB$. We define $f,F:\mathcal{NC}\lra\cD$ as follows:
 \begin{enumerate}
 \item[(a)] $f(\pi_1\cdot \pi_2)=f(\pi_1)\cdot f(\pi_2)$ and $F(\pi_1\cdot \pi_2)=F(\pi_1)\cdot F(\pi_2)$,\\
  where $\pi_1\cdot\pi_2\in NC(m+n)$ is obtained by juxtaposing $\pi_1\in NC(n)$ and $\pi_2\in NC(m)$.
     \item[(b)] $f(1_m)=\nu((\X b)^m)$ and $F(1_m)=\mu((\X b)^m)$,\\
      where $1_m\in NC(m)$ is the partition with a single block $(1, 2, \dots, m)$.
         \item[(c)] $f(\overline{|\pi_1|\pi_2|\dots|\pi_q|})=\nu(\X b \cdot f(\pi_1)\cdot \X b \cdot f(\pi_2)\cdots \X b \cdot f(\pi_q)\cdot \X b)$\\
         $F(\overline{|\pi_1|\pi_2|\dots|\pi_q|})=\mu(\X b \cdot f(\pi_1)\cdot \X b \cdot f(\pi_2)\cdots \X b \cdot f(\pi_q)\cdot \X b)$\\
         where $\overline{|\pi_1|\pi_2|\dots|\pi_q|}$ is the partition with a single exterior block with $q+1$ elements and with restrictions between the elements of the exterior block $\pi_1, \dots, \pi_q$.
 \end{enumerate}

 Similarly, we define $g,G:NC\lra\cD$ as above, replacing $\mu$ with ${}^c\rho_{\mu,\nu}$ and $\nu$ with $\rho_\nu$. With this notations, the moment-cumulant relation from  Definition \ref{def33} amounts to
 \begin{equation}\label{moebius1}
 F(\pi)=\sum_{\substack{\sigma\in \mathcal{NC}\\ \sigma\leq \pi}}G(\sigma)
 \end{equation}
 and applying the Moebius function property to equation \ref{moebius1}, we get
 \[G(\pi)=\sum_{\substack{\sigma\in \mathcal{NC}\\ \sigma\leq \pi}}F(\sigma)\cdot\text{moeb}
 (\sigma, \pi).
 \]

 Suppose now that $(\mu,\nu)\in\Sbd^0\times\Sb^0$ and that that $||\mu_n||_{cb}, ||\nu_n||_{cb}<M^{n+1}$ for all $n>0$, henceforth $||F(\sigma)||\leq M^{m+1}||b||^m$ for all $\sigma\in NC(m)$. Since ${}^c\kappa_{\mu,\nu,m}(b)=G(1_m)$, we have
 \begin{eqnarray*}
 \|{}^c\kappa_{\mu,\nu,m}(b)\|&\leq&
 \sum_{\sigma\in NC(m)}\|F(\sigma)\|\cdot|\text{moeb}
 (\sigma, 1_m)|\\
 &<&
 (\sharp NC(m))\cdot M^{m+1}\|b\|^m\cdot 4^m\\
 &<&
 \cdot M^{m+1}\|b\|^m\cdot 16^m,\ \text{since} \sharp NC(m)<4^m.
 \end{eqnarray*}
 Finally, $\kappa_{\nu, m}(b)={}^c\kappa_{\mu,\nu,m}(b)$, so the assertion is proven also for the $R$-transform.
\end{proof}

 \subsection{Main results}
 The following property of the $B$-transform is a non-commutative analogue of Proposition 3.1. from \cite{speicherwaroudi}, (i. e. the Nevalinna-Pick representation for the self-energy function of a real measure):
  \begin{thm}\label{selfenerg}
  Let $\cB$ be a C$^\ast$-algebra, $\cB\subset\cD$ be a unital inclusion of C$^\ast$-algebras and $\mu\in\Sbd$. Then there exists a selfadjoint $\alpha\in\cD$ and a $\mathbb{C}$-linear map $\sigma:\bx\lra\cD$, satisfying property (\ref{prop}) such that
    \[
    B_\mu(b)=\left[\alpha\cdot\mathds{1} +\widetilde{\sigma}\left(b(\mathds{1}-\X b)^{-1}\right)\right]\cdot b.
    \]
     Moreover, if $\mu\in\Sbd^0$, then the moments of $\nu$ do not grow faster than exponentially (i.e. $\nu$ satisfies property (\ref{pr2}).
  \end{thm}

   \begin{proof}
   Define the map $\beta_\mu:\bx\lra\cD$ as the $\cB$-bimodule extension of
   \[\beta_\mu(\X b_1 \X \cdots \X b_n)=B_{n,\mu}(b_1, \dots, b_n).\]

    Let $\cK=\bx\otimes_\cB\cD$. As shown in the proof of Theorem \ref{main1}, equation (\ref{bool2}), there exists some $\xi\in\cK$ and a selfadjoint map $T\in \cL(\cK)$ such that
    \[
    B_{n, \mu}(b_1, \dots, b_n)=\langle b_1T b_2\cdots T b_{n-1}(\xi b_n), \xi\rangle\]
    If $b_n=1$, the above relation becomes
    \begin{eqnarray*}
    \langle b_1T b_2\cdots T b_{n-1}\xi , \xi\rangle &=&B_{n, \mu}(b_1, \dots, b_{n-1}, 1)\\
    &=&\beta_{\mu}(\X b_1 \X \cdots \X b_{n-1}\X)
    \end{eqnarray*}
    Hence for all $f(\X)\in\bx$ we have that $\beta_\mu(\X f(\X)\X)=\langle f(T)\xi, \xi\rangle$. It follows that the map $\sigma:\bx\lra\cD$ given by $\sigma(f(\X))=\beta_\mu(\X f(\X)\X)$ satisfies property (\ref{prop}) and
    \[ B_{n, \mu}(b, \dots, b)=\sigma\left([\X b]^{n-1}\right)\cdot b,\]
    that is the conclusion for $\alpha=B_{1,\mu}$.
    The last part is a trivial consequence of the fact that $\mu\in\Sbd^0$ implies that $\beta_\mu\in\Sbd^0$.
   \end{proof}

   \begin{remark}\emph{If $\mu\in\Sbd^0$, the above result can be more explicit formulated in terms of the generalized rezolvent of $\mu$ from \cite{voi-qdif} and \cite{ncbp}. More precisely, for $\mu\in\Sbd^0$, its \emph{generalized rezolvent} or \emph{operator-valued Cauchy transform} is defined via}\end{remark}
   \[ \mathcal{G}_\mu(b)=\widetilde{\mu}\bigl([b-\mathds{1}\cdot\X]^{-1}\bigr).\]
 As shown in \cite{voi-qdif} and \cite{ncbp}, $\mathcal{G}_\mu$ is a non-commutative function, well-defined on $\mathbb{H}^{+}(\ncspace{\cB})$ and $\mathcal{G}_\mu(\mathbb{H}^{+}(M_n(\cB)))\subseteq \mathbb{H}^{-}(M_n(\cD))$. Hence, its reciprocal, $b\mapsto\bigl[\mathcal{G}_\mu(b)]^{-1}$ is also a non-commutative function, well defined on $\mathbb{H}^{+}(\ncspace{\cB})$. Moreover, identifying the coefficients, we have that
 \begin{equation}\label{rez-btransf}
 \mathds{1}-\bigl[\mathcal{G}_\mu(b)\bigr]^{-1}\cdot b^{-1}=B_\mu(b^{-1})
 \end{equation}
 for $b\in\mathbb{H}^{+}(\ncspace{B})$ with $||b^{-1}||<M$ for some $M>0$.

 Using equation \ref{rez-btransf} and Theorem \ref{selfenerg}, we obtain that, for $b$ as above,
 \begin{eqnarray*}
 b-\bigl[\mathcal{G}_\mu(b)\bigr]^{-1}
 &=&
 \alpha\cdot\mathds{1}+\widetilde{\sigma}(b^{-1}[\mathds{1}-\X\cdot b^{-1}]^{-1})\\
 &=&
 \alpha\cdot\mathds{1}+\widetilde{\sigma}([(\mathds{1}-\X \cdot b^{-1})\cdot b]^{-1})\\
 &=&
 \alpha\cdot\mathds{1}+\widetilde{\sigma}([b-\mathds{1}\cdot\X]^{-1}),
 \end{eqnarray*}
 The map $b\mapsto\widetilde{\sigma}([b-\mathds{1}\cdot\X]^{-1})$ extends to $\mathbb{H}^{+}(\ncspace{\cB})$ (see again \cite{voi-qdif} and \cite{voi09})  therefore we obtained that the \emph{operator-valued selfenergy function} of $\mu$ ($b\mapsto b-\bigl[\mathcal{G}_\mu(b)\bigr]^{-1}$) is the translate with a selfadjoint of the operator-valued Cauchy transform of some $\mathbb{C}$-linear map from $\bx$ to $\cD$ satisfying properties \ref{prop} and \ref{pr2}\hfill$\Box$

  For the main result of this section, Theorem \ref{repthm}, we will first need the following lemma:
  \begin{lemma}\label{auxlema} Let $\cB\subset\cD$ be a unital inclusion of C$^\ast$-algebras and $\rho:\bx\lra\cD$ be a unital $\cB$-bimodule map. Then the restriction of $\rho$ to $\bx_0$ satisfies property (\ref{prop}) if and only if there exists some ($\mathbb{C}$)-linear  map $\sigma=\sigma(\rho):\bx\lra\cD$ satisfying property (\ref{prop}) such that $\sigma(f(\X))=\rho(\X f(\X)\X)$ for all $f(\X)\in\bx$.
  \end{lemma}
  \begin{proof}
   Suppose that  first that the restriction of $\rho_{|\bx_0}$ satisfies (\ref{prop}) and define $\sigma$ via $\sigma(f(\X))=\rho(\X f(\X) \X)$.
   If  $\{ f_j(\X)\}_{j=1}^n$ is some family from $\bx_0$, then
   \begin{eqnarray*}
   \bigl[\sigma(f_j(\X)^\ast f_i(\X))\bigr]_{i,j=1}^n
   &=&\bigl[\rho(\X f_j(\X)^\ast f_i(\X)\X )\bigr]_{i,j=1}^n\\
   &=&\bigl[\rho\bigl([f_j(\X)\X]^\ast[f_i(\X)\X]\bigr)\bigl]_{i,j=1}^n\geq0
   \end{eqnarray*}

   For the converse, note that $\sigma$ satisfies (\ref{prop}) implies (cf \cite{paulsen}, pag. 42) that $\bx\otimes\cD$ is a semi-inner product $\cD$-module with respect to the pairing generated by
   \[\langle f(\X)\otimes d_1, g(\X)\otimes d_2\rangle =d_2^\ast \sigma(g(\X)^\ast f(\X))d_1.\]

    Fix now a family $\{f_j(\X)\}_{j=1}^n$ from $\bx_0$. Each $f_j(\X)$ can be written as
    \[ f_j(\X)=\sum_{k=1}^{N(j)}g_{k, j}(\X)\cdot\X\cdot \alpha_{k,j}\]
    with $g_{k,j}(\X)\in\bx$ and $\alpha_{k,j}\in\cB$.
    Denote $\eta_j=\sum_{k=1}^{N(j)}g_{k, j}(\X)\X \alpha_{k,j}\in\bx\otimes\cD$. Then
    \begin{eqnarray*}
    \bigr[ \rho(f_j(\X)^\ast f_i(\X))\bigl]_{i,j=1}^n
    &=&
    \bigr[\sum_{k=1}^{N(j)}\sum_{l=1}^{N(i)} \rho(\alpha_{k,j}^\ast\cdot \X\cdot g_{k,j}(\X)^\ast g_{l,i}(\X)\cdot\X\cdot\alpha_{l,i})\bigl]_{i,j=1}^n\\
    &=&
    \bigr[\sum_{k=1}^{N(j)}\sum_{l=1}^{N(i)}
    \alpha_{k,j}^\ast\cdot\sigma(g_{k,j}(\X)^\ast g_{l,i}(\X))\cdot\alpha_{l,i})\bigl]_{i,j=1}^n\\
    &=&
    \bigr[\langle \sum_{l=1}^{N(i)} g_{l,i}(\X)\otimes\alpha_{l,i}, \sum_{k=1}^{N(j)}g_{k,j}(\X)\otimes \alpha_{k,j}\rangle\bigl]_{i,j=1}^n\\
    &=&
    \bigr[\langle \eta_i, \eta_j\rangle\bigl]_{i,j=1}^n\geq 0.
    \end{eqnarray*}
    \end{proof}

An immediate consequence of   Theorem \ref{selfenerg} and the above Lemma is the following
   \begin{cor}
     If $\mu\in\Sbd$  then  $\beta_{\mu|\bx_0}$ satisfies property \emph{(\ref{prop})}.
   \end{cor}

 \begin{thm}\label{repthm}As before, $\cB\subset\cD$ will be a unital inclusion of unital C$\ast$-algebras.

 \emph{(i)}Let $\mu\in\Sb$. Then $\mu$ is $\boxplus$-infinitly divisible if and only if there exist some selfadjoint $\alpha\in\cB$ and some $\mathbb{C}$-linear map $\sigma:\bx\lra\cB$ satisfying property \emph{(\ref{prop})}, such that
  \[R_\mu(b)=\left[\alpha\cdot\mathds{1}+\widetilde{\sigma}\left(b(\mathds{1}-\X b)^{-1}\right)\right]\cdot b.\]

  \emph{(ii)} Let $(\mu,\nu)\in\Sbd\times\Sb$. Then  $(\mu,\nu)$ is $\cfree$-infinitely divisible if and only if $\nu$ is $\boxplus$-infinitely divisible and there exist some selfadjoint $\alpha\in\cB$ and some $\mathbb{C}$-linear map $\sigma:\bx\lra\cD$ satisfying property \emph{(\ref{prop})} such that
  \[{}^cR_{\mu,\nu}(b)=\left[\alpha\cdot\mathds{1}+\widetilde{\sigma}\left(b(\mathds{1}-\X b)^{-1}\right)\right]\cdot b.\]

  Moreover, if $\mu\in\Sb^0$, respectively $(\mu, \nu)\in\Sbd^0$, then the moments of the correspondent maps $\sigma$ do not grow faster than exponentially.
 \end{thm}
\begin{proof}

Since $\Sigma_{\cB:\cB}=\Sbd$, it suffices to prove the assertions for ${}^cR_{\mu,\nu}$.

 Since ${}^c\rho_{\mu,nu}(\X b_1\cdots \X b_n)=\kappa_{\mu,n}(b_1, \dots, b_n)$, Proposition \ref{tensorprop} implies
 that
 \[
 {}^cR_{\mu,\nu}(b)=M_{{}^c\rho_{\mu,\nu}}(b)-\mathds{1}.
 \]
  for all $b\in Nilp(\cB)$ or for $b$ or sufficiently small norm if $(\mu,\nu)\in\Sbd^0\times\Sb^0$.

  Identifying $\X$ to $\mathds{1}\cdot\X$, we have that
  \begin{eqnarray*}
 {}^cR_{\mu,\nu}(b)&=&\sum_{n=1}^\infty \widetilde{{}^c\rho_{\mu,\nu}}\bigl( (\X b)^n\bigr)\\
 &=&\widetilde{{}^c\rho_{\mu,\nu}}(\X)\cdot b+
 \left[\sum_{n=2}^\infty \widetilde{{}^c\rho_{\mu,\nu}}\bigl( (\X b)^{n-1}\X\bigr)\right]\cdot b\\
 &=&\left[\mu(\X)\cdot \mathds{1}
  +\sum_{n=2}^\infty \widetilde{{}^c\rho_{\mu,\nu}}\bigl(\X \cdot b(\X b)^{n-2}\cdot \X\bigr)
  \right]\cdot b
\end{eqnarray*}
 From Theorem \ref{thm45} and Lemma \ref{auxlema},  $(\mu,\nu)$ is $\cfree$-infinitely divisible if and only if there exists $\mathbb{C}$-linear  map $\sigma:\bx\lra\cD$ that satisfy property (\ref{prop}) such that $\sigma(f(\X))={}^c\rho_{\mu,\nu}(\X f(\X) \X)$ for all $f(\X)\in\bx$. That is
  \begin{eqnarray*}
  {}^cR_{\mu,\nu}(b)
  &=&
  \left[\mu(X)\cdot\mathds{1}+\sum_{l=0}^\infty\widetilde{\sigma}(b\cdot[\X b]^l)\right]\cdot b\\
  &=&
  \left[\mu(X)\cdot\mathds{1}+\widetilde{\sigma}(b[\mathds{1}-\X b]^{-1})\right]\cdot b
  \end{eqnarray*}
  hence the conclusion.

  For the last part, if $(\mu,\nu)\in\Sbd^0\times\Sb^0$, then Remark \ref{remark55} implies that $({}^c\rho_{\mu,\nu},\rho_\nu)$ is also in $\Sbd^0\times\Sb^0$, so they satisfy the condition (\ref{pr2}) for some $M>0$. The representation of ${}^cR_{\mu,\nu}$ gives (for $b_1,\dots, b_n\in M_m(\cB)$ and the identification $\X\cong\id_m\otimes\X$):

 \begin{eqnarray*}
 \|\widetilde{\sigma}(\X b_1 \X b_2\cdots b_n\X )\|&=&\|{}^c\rho_{\mu,\nu}(\X^2 b_1 \X b_2\cdots b_n\X^2)\|\\
 &<&M^{n+3}\|b_1\|\cdots\|b_n\|\\
 &<&(M^2+1)^{n+1}\|b_1\|\cdots\|b_n\|.
 \end{eqnarray*}

\end{proof}

%%%%%%%%%%%%%%%%%%%%%%%%%%%%%%%%%%%%%%%%%%%%%%%%%%%%%%%%%%%%%%%%%%%%%%%%%%%5
%%%%%%%%%%%%%%%%%%%%%%%%%%%%%%%%%%%%%%%%%%%%%%%%%%%%%%%%%%%%%%%%%%%%%%%%%%%%%%%%%%
%%%%%%%%%%%%%%%%%%%%%%%%%%%%%%%%%%%%%%%%%%%%%%%%%%%%%%%%%%%%%%%%%%%%
%%%%%%%%%%%%%%%%%%%%%%%%%%%%%%%%%%%%%%%%%%%%%%%%%%%%%%%%%%%%%%%%%%%%%%%%%%

\bibliographystyle{alpha}

\end{document}